\documentclass[a4paper,11pt]{amsart}

\usepackage{amsmath,amssymb,amsfonts}
\usepackage{amsthm}
\usepackage{geometry}
\usepackage{microtype}
\usepackage{aliascnt}
\usepackage[hidelinks]{hyperref}
\usepackage[nameinlink,capitalise]{cleveref}

\allowdisplaybreaks
\numberwithin{equation}{section}

\newtheorem{theorem}{Theorem}[section]
\newaliascnt{proposition}{theorem}
\newtheorem{proposition}[proposition]{Proposition}
\aliascntresetthe{proposition}
\newaliascnt{lemma}{theorem}
\newtheorem{lemma}[lemma]{Lemma}
\aliascntresetthe{lemma}
\newaliascnt{corollary}{theorem}

\aliascntresetthe{corollary}

\crefname{theorem}{Theorem}{Theorems}
\crefname{proposition}{Proposition}{Propositions}
\crefname{lemma}{Lemma}{Lemmas}
\crefname{corollary}{Corollary}{Corollaries}

\theoremstyle{definition}

\newaliascnt{example}{theorem} 
\newtheorem{example}[example]{Example} \aliascntresetthe{example}
\theoremstyle{remark}
\newtheorem{remark}[theorem]{Remark}

\crefname{example}{Example}{Examples} 
\Crefname{example}{Example}{Examples}

\newcommand{\N}{\mathbb N}
\newcommand{\cD}{\mathcal D}
\newcommand{\cP}{\mathcal P}
\newcommand{\cQ}{\mathcal Q}
\newcommand{\cR}{\mathcal R}
\newcommand{\cI}{\mathcal I}
\newcommand{\cV}{\mathcal V}
\newcommand{\ver}{\operatorname{ver}}
\newcommand{\transpose}{\mathsf T}
\title[Finite copula interpolation by linear programming]
{A linear programming approach to finite copula interpolation}

\author{Matja\v z Omladi\v c}
\address{Faculty of Mathematics and Physics, University of Ljubljana, and
Institute of Mathematics, Physics and Mechanics, Ljubljana, Slovenia}
\email{matjaz@omladic.net}

\author{Martin Vuk}
\address{Faculty of Computer and Information Science, University of Ljubljana,
Ljubljana, Slovenia}
\email{martin.vuk@fri.uni-lj.si}

\author{Alja\v z Zalar}
\address{Faculty of Computer and Information Science, University of Ljubljana;
Faculty of Mathematics and Physics, University of Ljubljana; and
Institute of Mathematics, Physics and Mechanics, Ljubljana, Slovenia}
\email{aljaz.zalar@fri.uni-lj.si}

\thanks{M.\ Omladi\v{c} acknowledges financial support from the
Slovenian Research and Innovation Agency through research core funding
No.\ P1-0448. A.\ Zalar acknowledges financial support from the
Slovenian Research and Innovation Agency through research core funding
No.\ P1-0288 and grants Nos.\ J1-50002, J1-60011, and J1-70017.}


\begin{document}

\begin{abstract}
We consider exact and interval-valued interpolation of multivariate
copulas at arbitrary finite configurations in the unit cube. We first
characterize when the prescribed data are compatible with a
quasi-copula and derive explicit extremal quasi-copula envelopes. The
copula interpolation problem is thereby reduced to a copula sandwich
problem between these envelopes. Restricting them to the canonical
rectangular grid generated by the data coordinates yields a finite
linear feasibility system whose elementary-volume constraints are
encoded by a cell--vertex incidence matrix. Farkas' lemma then identifies
the criterion of Omladi\v{c} and Stopar with an explicit certificate of
infeasibility. Although this criterion is indexed by infinitely many
formal cell combinations, we show that it suffices to test finitely many
vertices of a sign-polytopal subdivision. We also extend the construction
to \(k\)-increasing discrete quasi-copulas for \(2\leq k\leq d\).
\end{abstract}

\maketitle
\section{Introduction}

Let $E=\{\mathbf x_1,\ldots,\mathbf x_n\}\subseteq[0,1]^d$ be an
arbitrary finite set. The finite copula interpolation problem asks
whether there exists
a \(d\)-copula \(C\) such that
\begin{equation}\label{eq:intro-interpolation}
 C(\mathbf x_j)=q_j,
 \qquad
 j=1,\ldots,n,
\end{equation}
for prescribed values $q_1,\ldots,q_n\in[0,1]$.
More generally, given intervals $[a_j,b_j]\subseteq[0,1]$,
$j=1,\ldots,n$, the question is whether there exists a \(d\)-copula satisfying
\begin{equation}\label{eq:intro-interval}
 a_j\leq C(\mathbf x_j)\leq b_j,
 \qquad
 j=1,\ldots,n.
\end{equation}
No product-grid assumption is imposed on \(E\): the interpolation
points may form an arbitrary finite configuration in the unit cube.

Copula interpolation under finitely many pointwise constraints has
previously been studied in several forms. Best-possible bounds and
existence results for copulas with a prescribed value at one point are
closely related to the corresponding one-point quasi-copula bounds; see
\cite{RodriguezUbeda} and the references therein. The existence and
construction of multivariate copulas attaining prescribed quasi-copula
values at two arbitrary points were established in
\cite{KlementEtAl}. In the trivariate case, a linear-programming
approach was used in \cite{DeBaetsEtAl} to treat prescribed values at
several points: existence is guaranteed for two arbitrary points,
whereas it may fail for three points. These results motivate the
general finite interpolation problem considered here, in which both
the number and the coordinates of the interpolation points are
arbitrary.

The interval-valued problem \eqref{eq:intro-interval} is naturally
related to the more general problem of finding a function between
prescribed lower and upper bounds. The existence of a
\(k\)-increasing multivariate function satisfying $A\leq C\leq B$ has
been studied in the framework of coherence and avoidance of sure loss
for standardized functions and semicopulas; see
\cite{KlementEtAlCoherence} and the references therein. For copulas and
quasi-copulas, related discrete sandwich criteria were obtained by
Omladi\v{c} and Stopar in the bivariate setting
\cite{OmSt2020} and subsequently in
arbitrary dimension
\cite{OmSt2022}.

Our starting point is the compatibility of finite data with the class
of quasi-copulas. For exact data, the Fr\'echet--Hoeffding bounds and
the directed Lipschitz inequalities give necessary compatibility
conditions. We show that these conditions are also sufficient and
construct explicitly the smallest and the greatest quasi-copulas
interpolating the prescribed values. For interval-valued data,
represented by lower and upper functions
$A,B\colon E\to[0,1]$ satisfying $A\leq B$, we obtain analogous
extremal lower and upper quasi-copula envelopes, denoted by
\(\underline Q_A\) and \(\overline Q_B\). The original finite
interpolation constraints are then equivalent to the global sandwich
condition $\underline Q_A\leq Q\leq\overline Q_B$. Consequently, the
interval-valued copula interpolation problem is equivalent to the
existence of a copula between these two quasi-copula envelopes.

The finite set \(E\) need not itself have a rectangular structure.
Nevertheless, its coordinates generate a canonical rectangular grid $\cD_E$.
For each coordinate, we collect \(0\), \(1\), and all coordinate values
occurring among the points of \(E\), and then take the Cartesian
product of the resulting finite coordinate sets. Every point of \(E\)
is a node of the grid $\cD_E$. The generated grid is therefore an auxiliary
finite-dimensional representation of the original interpolation
problem, rather than an additional assumption on the geometry of the
data.
Restricting the extremal envelopes to the generated grid gives
discrete lower and upper bounds
$A_E:=\underline Q_A|_{\cD_E}$ and
$B_E:=\overline Q_B|_{\cD_E}$. The existence of a copula satisfying the
original interval constraints is thereby transformed into a finite
discrete copula sandwich problem. More precisely, we seek a discrete
copula \(C_E\) on \(\cD_E\) such that $A_E\leq C_E\leq B_E$.

We formulate this discrete sandwich problem in terms of values assigned
to the grid nodes. The elementary-cell volumes are encoded by a
cell--vertex incidence matrix \(H_E\). If
$\mathbf c=(C_E(\mathbf z_1),\ldots,C_E(\mathbf z_N))^\transpose$ is
the vector of node values, then $H_E\mathbf c\geq\mathbf0$ expresses
the nonnegativity of all elementary-cell volumes. The sandwich problem
is therefore equivalent to the finite feasibility system
$H_E\mathbf c\geq\mathbf0$,
$\mathbf a_E\leq\mathbf c\leq\mathbf b_E$, where
\(\mathbf a_E\) and \(\mathbf b_E\) are the node-value vectors of
\(A_E\) and \(B_E\), respectively.

Applying Farkas' lemma to this system gives a new interpretation of the
\(L^{(A,B)}\)-criterion introduced by Omladi\v{c} and Stopar
\cite{OmSt2020,OmSt2022}. A nonnegative multiplier vector assigns
weights to the elementary cells and determines a formal nonnegative
combination \(R\), whose node-multiplicity vector is obtained by
applying \(H_E^\transpose\) to the multiplier vector. The support
function of the interval box \([\mathbf a_E,\mathbf b_E]\) in this
direction is precisely \(L^{(A_E,B_E)}(R)\). We prove that the discrete
copula sandwich problem is feasible if and only if
\(L^{(A_E,B_E)}(R)\geq0\) for every formal nonnegative integer
combination \(R\) of elementary cells. Equivalently, a combination with
\(L^{(A_E,B_E)}(R)<0\) is an explicit certificate of infeasibility.
The canonical-grid construction therefore extends the discrete
sandwich criterion from bounds on a prescribed rectangular grid to
data on an arbitrary finite set.

Although the \(L^{(A_E,B_E)}\)-criterion is initially formulated
through infinitely many formal combinations, we show that it reduces
to finitely many inequalities. The associated dual functional
\(\Phi_E\) is continuous, positively homogeneous, and piecewise
linear. Positive homogeneity allows us to normalize every nonzero
nonnegative Farkas multiplier to the standard simplex. The signs of the
coordinates of the corresponding node-multiplicity vector subdivide
this simplex into finitely many sign polytopes. On each sign polytope,
\(\Phi_E\) is linear. It is therefore sufficient to verify its
nonnegativity at the vertices of all nonempty sign polytopes.

Finally, we extend this construction from \(d\)-increasing copulas to
\(k\)-increasing discrete quasi-copulas, where
\(k\in\{2,\ldots,d\}\). Elementary \(d\)-dimensional cells are replaced
by elementary \(k\)-dimensional coordinate boxes, whose volumes are
encoded by a \(k\)-cell--vertex incidence matrix \(H_E^{(k)}\). The
same Farkas argument yields a corresponding
\(L_k^{(A_E,B_E)}\)-criterion and a finite reduction to the vertices of
sign polytopes. For \(k=d\), this construction reduces to the discrete
copula sandwich criterion.

\subsection*{Organization of the paper}

\Cref{sec:preliminaries} fixes the notation and recalls the relevant
continuous and discrete notions. \Cref{sec:finite-data} establishes the
finite-data compatibility conditions and the extremal quasi-copula
envelopes. \Cref{sec:farkas} constructs the canonical grid, proves the
finite copula interpolation criterion, and illustrates its primal and
dual certificates. \Cref{sec:k-increasing-farkas} extends the argument
to \(k\)-increasing discrete quasi-copulas.

\section{Preliminaries}
\label{sec:preliminaries}

Throughout the paper, \(d\geq2\), inequalities between vectors are
understood coordinatewise, and
$[d]=\{1,\ldots,d\}.$

\subsection{Copulas and quasi-copulas}

For
\[
 \mathbf x=(x_1,\ldots,x_d),
 \qquad
 \mathbf y=(y_1,\ldots,y_d)\in[0,1]^d
\]
satisfying \(\mathbf x\leq\mathbf y\), a \(d\)-box is a set of the form
\[
 B=[\mathbf x,\mathbf y]
 :=
 \prod_{\ell=1}^d[x_\ell,y_\ell].
\]
The set of vertices of \(B\) is
\[
 \ver(B)
 :=
 \prod_{\ell=1}^d\{x_\ell,y_\ell\}.
\]
The box \(B\) is called \emph{nondegenerate} if \(x_\ell<y_\ell\) for every
\(\ell\in[d]\), and degenerate otherwise.

Let \(F\) be a real-valued function whose domain contains
\(\ver(B)\). If \(B\) is nondegenerate, the \(F\)-volume of \(B\) is
defined by
\begin{equation}\label{eq:F-volume}
 V_F(B)
 :=
 \sum_{\mathbf z\in\ver(B)}
 (-1)^{\#\{\ell\in[d]\colon z_\ell=x_\ell\}}F(\mathbf z).
\end{equation}
If \(B\) is degenerate, we put \(V_F(B):=0\). The function \(F\) is
called \emph{\(d\)-increasing} if
\[
 V_F(B)\geq0
\]
for every \(d\)-box \(B\subseteq[0,1]^d\).

A function \(C\colon[0,1]^d\to[0,1]\) is called a \emph{\(d\)-copula}, or
simply a \emph{copula}, if it satisfies the following conditions:
\begin{enumerate}
\item
\(C\) is \emph{grounded}, i.e.,
$C(u_1,\ldots,u_d)=0$
whenever \(u_\ell=0\) for some \(\ell\in[d]\);
\item
\(C\) has \emph{uniform one-dimensional marginals}, i.e.,
$
 C(1,\ldots,1,u_\ell,1,\ldots,1)=u_\ell
$
for every \(u_\ell\in[0,1]\) and every \(\ell\in[d]\);
\item
\(C\) is \(d\)-increasing.
\end{enumerate}
\smallskip

A function \(Q\colon[0,1]^d\to[0,1]\) is called a
\emph{\(d\)-quasi-copula}, or simply a \emph{quasi-copula}, if it satisfies the
following conditions:
\begin{enumerate}
\item
\(Q\) is grounded;
\item
\(Q\) has uniform one-dimensional margins;
\item
\(Q\) is increasing in each variable;
\item
\(Q\) is \emph{\(1\)-Lipschitz}, i.e.,
\[
 \lvert Q(\mathbf u)-Q(\mathbf v)\rvert
 \leq
 \sum_{\ell=1}^d\lvert u_\ell-v_\ell\rvert
\]
for all \(\mathbf u,\mathbf v\in[0,1]^d\).
\end{enumerate}
Here and throughout the paper, increasing means nondecreasing. Every
\(d\)-copula is a \(d\)-quasi-copula, but the converse does not hold
in general. For standard background, see
\cite{Nelsen,DuranteSempi}; for the lattice structure of multivariate
copulas and quasi-copulas, see \cite{FernandezNelsenUbeda}.
\smallskip

Every \(d\)-quasi-copula \(Q\) satisfies the
Fr\'echet--Hoeffding bounds
\cite[Section~3.3]{DuranteSempi}; see also \cite{AriasEtAl}:
\[
 W_d(\mathbf u)
 \leq Q(\mathbf u)
 \leq M_d(\mathbf u),
 \qquad
 \mathbf u\in[0,1]^d,
\]
where
\begin{equation}
\label{eq:FH-bounds}
 W_d(\mathbf u)
 :=
 \max\left\{
  0,\sum_{\ell=1}^d u_\ell-d+1
 \right\},
 \qquad
 M_d(\mathbf u)
 :=
 \min_{1\leq\ell\leq d}u_\ell.
\end{equation}

\subsection{\texorpdfstring{\(k\)-increasing functions}{k-increasing functions}}

We will also use lower-dimensional coordinate boxes. Let
\(J\subseteq[d]\) be a nonempty set of indices. A
\(|J|\)-dimensional coordinate box is
a set of the form
\[
 B_J
 =
 \prod_{\ell=1}^d I_\ell,
\]
where
\[
 I_\ell
 =
 \begin{cases}
 [x_\ell,y_\ell], & \ell\in J,\\
 \{q_\ell\},      & \ell\notin J,
 \end{cases}
\]
for some \(x_\ell<y_\ell\) when \(\ell\in J\) and some fixed
\(q_\ell\in[0,1]\) when \(\ell\notin J\). The coordinates indexed by
\(J\) are called \emph{active}, and the dimension of \(B_J\) is \(|J|\).

The set of vertices of \(B_J\) is
\[
 \ver(B_J)
 =
 \left\{
  \mathbf z\in[0,1]^d:
  z_\ell\in\{x_\ell,y_\ell\}\text{ for }\ell\in J,
  \quad
  z_\ell=q_\ell\text{ for }\ell\notin J
 \right\}.
\]
For a real-valued function \(F\) whose domain contains
\(\ver(B_J)\), the \(F\)-volume of \(B_J\) is defined by
\begin{equation}\label{eq:F-volume-coordinate-box}
 V_F(B_J)
 :=
 \sum_{\mathbf z\in\ver(B_J)}
 (-1)^{\#\{\ell\in J:z_\ell=x_\ell\}}F(\mathbf z).
\end{equation}

For \(k\in[d]\), the function \(F\) is called \emph{\(k\)-increasing} if
\[
 V_F(B_J)\geq0
\]
for every \(k\)-dimensional coordinate box \(B_J\) with \(|J|=k\).

\subsection{Discrete copulas and quasi-copulas}

Let
\[
 \cD=D_1\times\cdots\times D_d\subseteq[0,1]^d
\]
be a finite rectangular grid, where \(0,1\in D_\ell\) for every
\(\ell\in[d]\). A function \(Q\colon\cD\to[0,1]\) is called a
\emph{discrete \(d\)-quasi-copula} if it satisfies the following conditions:
\begin{enumerate}
\item
\(Q\) is grounded, i.e., \(Q(\mathbf u)=0\) whenever some coordinate
of \(\mathbf u\) is zero;
\item
\(Q\) has uniform one-dimensional margins, i.e.,
$
 Q(1,\ldots,1,u_\ell,1,\ldots,1)=u_\ell
$
for every \(u_\ell\in D_\ell\) and every \(\ell\in[d]\);
\item
\(Q\) is increasing in each variable on \(\cD\);
\item
\(Q\) satisfies
\[
 \lvert Q(\mathbf u)-Q(\mathbf v)\rvert
 \leq
 \sum_{\ell=1}^d\lvert u_\ell-v_\ell\rvert
\]
for all \(\mathbf u,\mathbf v\in\cD\).
\end{enumerate}

Equivalently, the last two conditions may be checked along coordinate
edges: whenever \(\mathbf u,\mathbf v\in\cD\) differ only in
coordinate \(\ell\) and \(u_\ell\leq v_\ell\), one has
\[
 0
 \leq
 Q(\mathbf v)-Q(\mathbf u)
 \leq
 v_\ell-u_\ell.
\]

For a coordinate box \(B_J\) whose vertices belong to \(\cD\), its
\(Q\)-volume is defined by
\eqref{eq:F-volume-coordinate-box}. A discrete quasi-copula \(Q\) is
called \(k\)-increasing if
\[
 V_Q(B_J)\geq0
\]
for every \(k\)-dimensional coordinate box \(B_J\) with
\(|J|=k\) and \(\ver(B_J)\subseteq\cD\). A \emph{discrete \(d\)-copula} is a
\(d\)-increasing discrete \(d\)-quasi-copula.

\begin{proposition}[Extension of a discrete copula]
\label{prop:discrete-copula-extension}
Let \(\cD=D_1\times\cdots\times D_d\subseteq[0,1]^d\) be a finite
rectangular grid such that \(0,1\in D_\ell\) for every \(\ell\in[d]\).
Every discrete \(d\)-copula
\(C_{\cD}\colon\cD\to[0,1]\) admits a \(d\)-copula extension
\(C\colon[0,1]^d\to[0,1]\) satisfying \(C|_{\cD}=C_{\cD}\).
\end{proposition}

\begin{proof}
A discrete \(d\)-copula on \(\cD\) is a \(d\)-subcopula whose domain
is the product set \(D_1\times\cdots\times D_d\). The conclusion follows
from the extension theorem for subcopulas; see
\cite[Section~2.3]{DuranteSempi}.
\end{proof}

\section{Finite interpolation data and quasi-copula envelopes}
\label{sec:finite-data}

We first characterize finite exact data that admit a quasi-copula
interpolant and construct the least and greatest such interpolants.
We then derive analogous envelopes for interval-valued data. These
envelopes turn the original pointwise constraints into the global
sandwich bounds used in \cref{sec:farkas}.

For \(\mathbf u,\mathbf v\in[0,1]^d\), define the directed distance
from \(\mathbf v\) to \(\mathbf u\) by
\begin{equation}\label{eq:directed-distance}
 d_+(\mathbf u,\mathbf v)
 := \sum_{\ell=1}^d (u_\ell-v_\ell)_+,
 \qquad
 t_+:=\max\{t,0\}.
\end{equation}
Recall the Fr\'echet--Hoeffding bounds
\(W_d\) and \(M_d\) from \eqref{eq:FH-bounds}.

\subsection{Exact data}

Let $q_j\in[0,1]$ be prescribed at $\mathbf x_j$. We call the exact
data $(\mathbf x_j,q_j)_{j=1}^n$ \emph{quasi-copula compatible} if
\begin{align}
 W_d(\mathbf x_j)&\leq q_j\leq M_d(\mathbf x_j),
 &&j=1,\ldots,n,\label{eq:exact-FH}\\
 q_i-q_j&\leq d_+(\mathbf x_i,\mathbf x_j),
 &&i,j=1,\ldots,n.\label{eq:exact-directed}
\end{align}
Define
\begin{align}
 \underline Q(\mathbf u)
 &:=\max\left\{
 W_d(\mathbf u),
 \max_{1\leq i\leq n}
 \bigl[q_i-d_+(\mathbf x_i,\mathbf u)\bigr]
 \right\},\label{eq:exact-lower}\\
 \overline Q(\mathbf u)
 &:=\min\left\{
 M_d(\mathbf u),
 \min_{1\leq i\leq n}
 \bigl[q_i+d_+(\mathbf u,\mathbf x_i)\bigr]
 \right\}.\label{eq:exact-upper}
\end{align}

The one-point versions of the bounds below are due to
Rodr\'iguez-Lallena and \'Ubeda-Flores
\cite[Theorem~3.2]{RodriguezUbeda}; see also
\cite[Theorem~2.1]{KlementEtAl}. The finite-point result follows by
combining those bounds with the fact that finite pointwise maxima and
minima of \(d\)-quasi-copulas remain \(d\)-quasi-copulas. We include
the details because the extremal envelopes are used below.

\begin{proposition}[Finite-point quasi-copula interpolation]
\label{prop:exact-envelopes}
The exact data admit a quasi-copula interpolant if and only if
\eqref{eq:exact-FH}--\eqref{eq:exact-directed} hold. In that case,
\(\underline Q\) and \(\overline Q\) are quasi-copulas,
\[
 \underline Q(\mathbf x_j)=q_j=\overline Q(\mathbf x_j),
 \qquad j=1,\ldots,n,
\]
and every quasi-copula \(Q\) satisfying \(Q(\mathbf x_j)=q_j\) obeys
\[
 \underline Q\leq Q\leq\overline Q.
\]
\end{proposition}

\begin{proof}
Assume first that a quasi-copula \(Q\) interpolating the prescribed data
exists. Every \(d\)-quasi-copula satisfies the Fr\'echet--Hoeffding
inequalities
\[
W_d(\mathbf u)
\leq Q(\mathbf u)
\leq M_d(\mathbf u),
\qquad
\mathbf u\in[0,1]^d.
\]
Evaluating these inequalities at \(\mathbf u=\mathbf x_j\) gives
\eqref{eq:exact-FH}.

Moreover, coordinatewise monotonicity and the Lipschitz property of a
quasi-copula imply the directed inequality
\begin{equation}
\label{eq:directed-Lipschitz-proof}
Q(\mathbf u)-Q(\mathbf v)
\leq
\sum_{\ell=1}^d(u_\ell-v_\ell)_+
=
d_+(\mathbf u,\mathbf v),
\qquad
\mathbf u,\mathbf v\in[0,1]^d.
\end{equation}
Indeed, let
\[
 \mathbf w:=\mathbf u\vee\mathbf v
 =\bigl(\max\{u_1,v_1\},\ldots,\max\{u_d,v_d\}\bigr).
\]
coordinatewise monotonicity and the Lipschitz condition give
\[
Q(\mathbf u)-Q(\mathbf v)
\leq
Q(\mathbf w)-Q(\mathbf v)
\leq
\sum_{\ell=1}^d(w_\ell-v_\ell)
=
\sum_{\ell=1}^d(u_\ell-v_\ell)_+.
\]
Applying \eqref{eq:directed-Lipschitz-proof} to
$\mathbf u=\mathbf x_i$ and $\mathbf v=\mathbf x_j$ yields
\eqref{eq:exact-directed}.

Conversely, assume that
\eqref{eq:exact-FH}--\eqref{eq:exact-directed} hold. For every
$i\in\{1,\ldots,n\}$, define
\begin{align*}
Q_i^-(\mathbf u)
&:=
\max\left\{
W_d(\mathbf u),
q_i-d_+(\mathbf x_i,\mathbf u)
\right\},\\
Q_i^+(\mathbf u)
&:=
\min\left\{
M_d(\mathbf u),
q_i+d_+(\mathbf u,\mathbf x_i)
\right\}.
\end{align*}
Since
$
W_d(\mathbf x_i)
\leq q_i
\leq M_d(\mathbf x_i),
$
the one-point quasi-copula bounds \cite[Theorem~3.2]{RodriguezUbeda}
imply that $Q_i^-$ and $Q_i^+$ are quasi-copulas,
that
$
Q_i^-(\mathbf x_i)
=
q_i
=
Q_i^+(\mathbf x_i),
$  
and that every quasi-copula $Q$ satisfying
$Q(\mathbf x_i)=q_i$ fulfills
$
    Q_i^-\leq Q\leq Q_i^+.
$

Finite pointwise maxima and minima of \(d\)-quasi-copulas are again
\(d\)-quasi-copulas. Indeed, groundedness, uniform one-dimensional
margins, coordinatewise monotonicity, and the common \(1\)-Lipschitz
bound are preserved under both operations. Consequently,
\[
 \underline Q = \max_{1\leq i\leq n}Q_i^-,
 \qquad
 \overline Q = \min_{1\leq i\leq n}Q_i^+,
\]
are quasi-copulas. These are precisely the functions defined in
\eqref{eq:exact-lower} and \eqref{eq:exact-upper}.

It remains to verify that they interpolate all prescribed values.
For $i,j\in\{1,\ldots,n\}$, condition
\eqref{eq:exact-directed} gives
$
q_i-d_+(\mathbf x_i,\mathbf x_j)
\leq q_j.
$
Together with
$
W_d(\mathbf x_j)\leq q_j,
$
this yields
$
Q_i^-(\mathbf x_j)\leq q_j.
$
For $i=j$, equality holds:
$
Q_j^-(\mathbf x_j)=q_j.
$
Therefore,
\[
\underline Q(\mathbf x_j)
=
\max_{1\leq i\leq n}Q_i^-(\mathbf x_j)
=
q_j.
\]

Similarly, applying \eqref{eq:exact-directed} to the ordered pair
$(j,i)$ gives
$
q_j-q_i
\leq
d_+(\mathbf x_j,\mathbf x_i),
$
or equivalently,
$
q_j
\leq
q_i+d_+(\mathbf x_j,\mathbf x_i).
$
Since
$
q_j\leq M_d(\mathbf x_j),
$
we obtain
$
q_j\leq Q_i^+(\mathbf x_j).
$
Again, equality holds for $i=j$, and hence
\[
\overline Q(\mathbf x_j)
=
\min_{1\leq i\leq n}Q_i^+(\mathbf x_j)
=
q_j.
\]

Finally, let $Q$ be any quasi-copula interpolating the prescribed
data. Applying the one-point extremal bounds at every
$(\mathbf x_i,q_i)$ gives
$Q_i^-\leq Q\leq Q_i^+$,
\(i=1,\ldots,n\).
Taking the pointwise maximum of the lower bounds and the pointwise
minimum of the upper bounds yields
\[
\underline Q
=
\max_{1\leq i\leq n}Q_i^-
\leq Q
\leq
\min_{1\leq i\leq n}Q_i^+
=
\overline Q.
\]
Thus, \(\underline Q\) and \(\overline Q\) are, respectively, the least
and greatest quasi-copulas interpolating the prescribed data in the
pointwise order.
\end{proof}

\subsection{Interval-valued data}

Let \(0\leq a_j\leq b_j\leq1\). We call
$(\mathbf x_j;[a_j,b_j])_{j=1}^n$ \emph{quasi-copula compatible} if
\begin{align}
 a_i&\leq M_d(\mathbf x_i),&&i=1,\ldots,n,\label{eq:int-upper-FH}\\
 W_d(\mathbf x_j)&\leq b_j,&&j=1,\ldots,n,\label{eq:int-lower-FH}\\
 a_i-b_j&\leq d_+(\mathbf x_i,\mathbf x_j),&&i,j=1,\ldots,n.
 \label{eq:int-directed}
\end{align}
Define the interval envelopes
\begin{align}
 \underline Q_A(\mathbf u)
 &=\max\left\{
 W_d(\mathbf u),
 \max_{1\leq i\leq n}
 \bigl[a_i-d_+(\mathbf x_i,\mathbf u)\bigr]
 \right\},\label{eq:interval-lower}\\
 \overline Q_B(\mathbf u)
 &=\min\left\{
 M_d(\mathbf u),
 \min_{1\leq j\leq n}
 \bigl[b_j+d_+(\mathbf u,\mathbf x_j)\bigr]
 \right\}.\label{eq:interval-upper}
\end{align}

\begin{proposition}[Extremal interval envelopes]
\label{prop:interval-envelopes}
The interval data admit a quasi-copula $Q$ satisfying
\[
 a_j\leq Q(\mathbf x_j)\leq b_j,
 \qquad j=1,\ldots,n,
\]
if and only if \eqref{eq:int-upper-FH}--\eqref{eq:int-directed} hold.
In that case, \(\underline Q_A\) and \(\overline Q_B\) are quasi-copulas,
\(\underline Q_A\leq\overline Q_B\), and a quasi-copula \(Q\) satisfies
the interval constraints if and only
if
\begin{equation}\label{eq:interval-global-sandwich}
 \underline Q_A\leq Q\leq\overline Q_B.
\end{equation}
\end{proposition}

\begin{proof}
If \(Q\) satisfies the interval constraints, the Fr\'echet--Hoeffding and
directed Lipschitz inequalities give
\eqref{eq:int-upper-FH}--\eqref{eq:int-directed}. They also give, for
every \(\mathbf u\),
\[
 Q(\mathbf u)\geq a_i-d_+(\mathbf x_i,\mathbf u)
\qquad
\text{and}
\qquad
 Q(\mathbf u)\leq b_j+d_+(\mathbf u,\mathbf x_j),
\]
so \eqref{eq:interval-global-sandwich} follows.

Conversely, suppose that the compatibility conditions hold. Put
\[
 \widetilde a_i:=\max\{a_i,W_d(\mathbf x_i)\},
 \qquad
 \widetilde b_i:=\min\{b_i,M_d(\mathbf x_i)\}.
\]
Then
\[
 W_d(\mathbf x_i)\leq\widetilde a_i\leq M_d(\mathbf x_i),
 \qquad
 W_d(\mathbf x_i)\leq\widetilde b_i\leq M_d(\mathbf x_i).
\]
If \(a_i<W_d(\mathbf x_i)\), the directed Lipschitz inequality for
\(W_d\) gives
\[
 a_i-d_+(\mathbf x_i,\mathbf u)
 \leq W_d(\mathbf x_i)-d_+(\mathbf x_i,\mathbf u)
 \leq W_d(\mathbf u),
\]
so both the original and clipped terms in \eqref{eq:interval-lower}
are redundant.
Similarly, if \(b_i>M_d(\mathbf x_i)\), then
\[
 M_d(\mathbf u)
 \leq M_d(\mathbf x_i)+d_+(\mathbf u,\mathbf x_i)
 \leq b_i+d_+(\mathbf u,\mathbf x_i),
\]
so both the original and clipped terms in \eqref{eq:interval-upper}
are redundant.
Thus the envelopes are unchanged if \(a_i,b_i\) are replaced by
\(\widetilde a_i,\widetilde b_i\), respectively. The one-point bounds
\cite[Theorem~3.2]{RodriguezUbeda} therefore show that every function
entering the two finite extrema is a quasi-copula. Since finite
pointwise maxima and minima preserve the quasi-copula properties,
\(\underline Q_A\) and \(\overline Q_B\) are quasi-copulas.

It remains to verify that
\(\underline Q_A\leq\overline Q_B\).
Fix \(\mathbf u\in[0,1]^d\). First,
\(W_d(\mathbf u)\leq M_d(\mathbf u)\). Moreover, the directed
Lipschitz inequalities for the Fr\'echet--Hoeffding bounds give
\[
 W_d(\mathbf u)
 \leq W_d(\mathbf x_j)+d_+(\mathbf u,\mathbf x_j)
 \leq b_j+d_+(\mathbf u,\mathbf x_j)
\]
for every \(j\), and
\[
 a_i-d_+(\mathbf x_i,\mathbf u)
 \leq M_d(\mathbf x_i)-d_+(\mathbf x_i,\mathbf u)
 \leq M_d(\mathbf u)
\]
for every \(i\). Finally, by \eqref{eq:int-directed} and the directed
triangle inequality,
\[
\begin{aligned}
 a_i-d_+(\mathbf x_i,\mathbf u)
 &\leq
 b_j+d_+(\mathbf x_i,\mathbf x_j)
      -d_+(\mathbf x_i,\mathbf u)
    \leq
 b_j+d_+(\mathbf u,\mathbf x_j)
\end{aligned}
\]
for every \(i,j\). Thus every term in the maximum defining
\(\underline Q_A(\mathbf u)\) is bounded above by every term in the
minimum defining \(\overline Q_B(\mathbf u)\). Hence
\(\underline Q_A(\mathbf u)\leq\overline Q_B(\mathbf u)\).
Since \(\mathbf u\) was arbitrary, this proves
\(\underline Q_A\leq\overline Q_B\).

At every interpolation point \(\mathbf x_j\), the terms with index
\(j\) in \eqref{eq:interval-lower} and \eqref{eq:interval-upper} give
\(\underline Q_A(\mathbf x_j)\geq a_j\) and
\(\overline Q_B(\mathbf x_j)\leq b_j\).
Consequently, every quasi-copula \(Q\) satisfying
\(\underline Q_A\leq Q\leq\overline Q_B\)
also satisfies
\(a_j\leq Q(\mathbf x_j)\leq b_j\),
\(j=1,\ldots,n\).
In particular, both \(\underline Q_A\) and \(\overline Q_B\) satisfy the
prescribed interval constraints.
\end{proof}

\begin{remark}\label{rem:copula-sandwich}
Under the compatibility assumptions of \cref{prop:interval-envelopes},
the finite interval interpolation problem is feasible if and only if
there exists a \(d\)-copula \(C\) such that
\[
 \underline Q_A\leq C\leq\overline Q_B.
\]
For exact data, the analogous bounds are \(\underline Q\) and
\(\overline Q\), which coincide at every interpolation point.
\end{remark}

\section{Copula sandwich problems and Farkas duality}
\label{sec:farkas}

We reformulate finite copula interpolation as a linear feasibility
problem in the node values of a canonical rectangular grid.
Linear-programming formulations for prescribed copula values were used
in the trivariate setting in \cite{DeBaetsEtAl}; discrete sandwich
problems and the associated \(L^{(A,B)}\)-criterion were studied in
\cite{OmSt2020,OmSt2022}. Here a cell--vertex incidence matrix encodes
all elementary-cell volumes. Farkas' lemma identifies a negative
\(L\)-value with an infeasibility certificate, while a sign-polytopal
subdivision reduces the resulting infinite family of inequalities to a
finite vertex test. The final examples give complementary primal and
dual certificates.

Let
\[
 E=\{\mathbf x_1,\ldots,\mathbf x_n\}\subseteq[0,1]^d,
 \qquad
 \mathbf x_j=(x_{1j},\ldots,x_{dj}).
\]

The finite set \(E\) need not itself form a rectangular grid. We therefore
begin by constructing the smallest rectangular grid containing it. For each
coordinate \(\ell\in\{1,\ldots,d\}\), define
\begin{equation}\label{eq:coordinate-set}
 D_\ell(E)=\{0,1,x_{\ell 1},\ldots,x_{\ell n}\},
\end{equation}
where repeated values are omitted. Enumerate the resulting set in increasing
order as
\[
 0=t_{\ell,0}<t_{\ell,1}<\cdots<t_{\ell,n_\ell}=1.
\]
The rectangular grid generated by \(E\) is
\begin{equation}\label{eq:generated-grid}
 \cD_E=D_1(E)\times\cdots\times D_d(E).
\end{equation}
By construction, \(E\subseteq\cD_E\), even though \(E\) itself need not have
a product structure.

The elementary cells of \(\cD_E\) are indexed by
\[
 \cI_E=\prod_{\ell=1}^d\{1,\ldots,n_\ell\}.
\]
For each multi-index
\(\mathbf i=(i_1,\ldots,i_d)\in\cI_E\), let
\begin{equation}\label{eq:elementary-cell}
 B_{\mathbf i}
 =
 \prod_{\ell=1}^d
 [t_{\ell,i_{\ell}-1},t_{\ell,i_\ell}].
\end{equation}

This formulation connects the \(L^{(A,B)}\)-criterion of Omladi\v{c}
and Stopar \cite[Proposition~16]{OmSt2020} and its multivariate
extension \cite[Proposition~14]{OmSt2022} to standard linear-programming
duality.

Enumerate the nodes of \(\cD_E\) as
\[
 \cD_E=\{\mathbf z_1,\ldots,\mathbf z_N\}.
\]
For \(\mathbf i\in\cI_E\), write
\[
 B_{\mathbf i}
 =
 [\mathbf x^{\mathbf i},\mathbf y^{\mathbf i}]
 =
 \prod_{\ell=1}^d
 [t_{\ell,i_{\ell}-1},t_{\ell,i_\ell}],
\qquad\text{where}
\quad
 x^{\mathbf i}_\ell=t_{\ell,i_{\ell}-1},
 \;
 y^{\mathbf i}_\ell=t_{\ell,i_\ell}.
\]
Define the vertex-multiplicity function of \(B_{\mathbf i}\) by
\begin{equation}\label{eq:multiplicity}
 m_{\mathbf i}(\mathbf z)
 =
 \begin{cases}
 (-1)^{\#\{\ell:z_\ell=x^{\mathbf i}_\ell\}},
     & \mathbf z\in\ver(B_{\mathbf i}),\\
 0,  & \mathbf z\notin\ver(B_{\mathbf i}).
 \end{cases}
\end{equation}
The corresponding cell--vertex incidence matrix is
\begin{equation}\label{eq:H}
 H_E
 =
 \bigl[m_{\mathbf i}(\mathbf z_s)\bigr]_{
   \mathbf i\in\cI_E,\;1\leq s\leq N}
 \in\{-1,0,1\}^{|\cI_E|\times N}.
\end{equation}

If
\[
 \mathbf c
 =
 \bigl(c(\mathbf z_1),\ldots,c(\mathbf z_N)\bigr)^\transpose
\]
is a vector of values assigned to the grid nodes, then
\[
 (H_E\mathbf c)_{\mathbf i}
 =
 \sum_{\mathbf z\in\ver(B_{\mathbf i})}
 m_{\mathbf i}(\mathbf z)c(\mathbf z)
\]
is the volume assigned to the elementary cell \(B_{\mathbf i}\).
Consequently,
$
 H_E\mathbf c\geq0
$
is equivalent to requiring that every elementary cell have nonnegative
volume. Since every box with vertices in \(\cD_E\) is a union of elementary
cells, finite additivity implies that every such box also has nonnegative
volume.

\begin{example}[A bivariate grid generated by two incomparable points]
\label{ex:bivariate-incidence}
Let
\[
 E=
 \left\{
 \left(\frac14,\frac34\right),
 \left(\frac23,\frac13\right)
 \right\}
 \subseteq(0,1)^2.
\]
Note that the two points are incomparable in the coordinatewise order.
They are also nonsymmetric in the sense that \(E\) is not invariant
under interchange of the two coordinates.

The coordinate sets generated by \(E\) are
\[
 D_1(E)
 =
 \left\{0,\frac14,\frac23,1\right\},
 \qquad
 D_2(E)
 =
 \left\{0,\frac13,\frac34,1\right\}.
\]
Hence
\[
 \cD_E
 =
 \left\{0,\frac14,\frac23,1\right\}
 \times
 \left\{0,\frac13,\frac34,1\right\}.
\]

Order the sixteen grid nodes lexicographically:
\[
\begin{aligned}
 \mathbf z_1&=(0,0),&
 \mathbf z_2&=\left(0,\frac13\right),&
 \mathbf z_3&=\left(0,\frac34\right),&
 \mathbf z_4&=(0,1),\\
 \mathbf z_5&=\left(\frac14,0\right),&
 \mathbf z_6&=\left(\frac14,\frac13\right),&
 \mathbf z_7&=\left(\frac14,\frac34\right),&
 \mathbf z_8&=\left(\frac14,1\right),\\
 \mathbf z_9&=\left(\frac23,0\right),&
 \mathbf z_{10}&=\left(\frac23,\frac13\right),&
 \mathbf z_{11}&=\left(\frac23,\frac34\right),&
 \mathbf z_{12}&=\left(\frac23,1\right),\\
 \mathbf z_{13}&=(1,0),&
 \mathbf z_{14}&=\left(1,\frac13\right),&
 \mathbf z_{15}&=\left(1,\frac34\right),&
 \mathbf z_{16}&=(1,1).
\end{aligned}
\]

For \(i,j\in\{1,2,3\}\), let
\[
 B_{ij}
 =
 [t_{1,i-1},t_{1,i}]
 \times
 [t_{2,j-1},t_{2,j}],
\]
where
\[
 (t_{1,0},t_{1,1},t_{1,2},t_{1,3})
 =
 \left(0,\frac14,\frac23,1\right)
\qquad
\text{and}
\qquad
 (t_{2,0},t_{2,1},t_{2,2},t_{2,3})
 =
 \left(0,\frac13,\frac34,1\right).
\]
Thus, the nine elementary cells are ordered as
\[
 B_{11},B_{12},B_{13},
 B_{21},B_{22},B_{23},
 B_{31},B_{32},B_{33}.
\]

With this ordering, the cell--vertex incidence matrix is
\[
 H_E
 =
 \begin{pmatrix}
  1&-1& 0& 0&-1& 1& 0& 0& 0& 0& 0& 0& 0& 0& 0& 0\\
  0& 1&-1& 0& 0&-1& 1& 0& 0& 0& 0& 0& 0& 0& 0& 0\\
  0& 0& 1&-1& 0& 0&-1& 1& 0& 0& 0& 0& 0& 0& 0& 0\\
  0& 0& 0& 0& 1&-1& 0& 0&-1& 1& 0& 0& 0& 0& 0& 0\\
  0& 0& 0& 0& 0& 1&-1& 0& 0&-1& 1& 0& 0& 0& 0& 0\\
  0& 0& 0& 0& 0& 0& 1&-1& 0& 0&-1& 1& 0& 0& 0& 0\\
  0& 0& 0& 0& 0& 0& 0& 0& 1&-1& 0& 0&-1& 1& 0& 0\\
  0& 0& 0& 0& 0& 0& 0& 0& 0& 1&-1& 0& 0&-1& 1& 0\\
  0& 0& 0& 0& 0& 0& 0& 0& 0& 0& 1&-1& 0& 0&-1& 1
 \end{pmatrix}.
\]
For the node-value vector
\[
 \mathbf c
 =
 \bigl(c(\mathbf z_1),\ldots,c(\mathbf z_{16})\bigr)^\transpose,
\]
we therefore have
\[
 H_E\mathbf c
 =
 \begin{pmatrix}
 V_c(B_{11})\\
 V_c(B_{12})\\
 V_c(B_{13})\\
 V_c(B_{21})\\
 V_c(B_{22})\\
 V_c(B_{23})\\
 V_c(B_{31})\\
 V_c(B_{32})\\
 V_c(B_{33})
 \end{pmatrix}.
\]
Consequently,
\[
 H_E\mathbf c\geq\mathbf0_9
\]
is precisely the requirement that all nine elementary cells have
nonnegative \(c\)-volume.\hfill $\triangle$
\end{example}

We next consider formal nonnegative integer combinations of elementary
cells. Let
\begin{equation}\label{eq:formal-family}
 \cR_E
 =
 \left\{
  \bigsqcup_{\mathbf i\in\cI_E}
  \lambda_{\mathbf i}B_{\mathbf i}
  :
  \boldsymbol\lambda
  =
  (\lambda_{\mathbf i})_{\mathbf i\in\cI_E}
  \in\N_0^{\cI_E}
 \right\}.
\end{equation}
Here, \(\bigsqcup\) denotes a formal multiset union: the coefficient
\(\lambda_{\mathbf i}\) records the multiplicity with which the elementary
cell \(B_{\mathbf i}\) occurs. Thus, the elements of \(\cR_E\) need not be
geometric boxes or disjoint unions of cells.

For
\begin{equation}
\label{def:R}
 R\equiv R_\lambda
 :=
 \bigsqcup_{\mathbf i\in\cI_E}
 \lambda_{\mathbf i}B_{\mathbf i}
 \in\cR_E,
\end{equation}
define its node-multiplicity function by
\begin{equation}\label{eq:formal-multiplicity}
 m_R(\mathbf z)
 =
 \sum_{\mathbf i\in\cI_E}
 \lambda_{\mathbf i}m_{\mathbf i}(\mathbf z),
 \qquad
 \mathbf z\in\cD_E.
\end{equation}
The corresponding multiplicity vector is
\[
 \mathbf m_R
 =
 \bigl(m_R(\mathbf z_1),\ldots,m_R(\mathbf z_N)\bigr)^\transpose.
\]
After fixing an ordering of \(\cI_E\) and regarding
\(\boldsymbol\lambda\) as a column vector in
\(\N_0^{|\cI_E|}\), we have
\begin{equation}\label{eq:m-vector}
 \mathbf m_R=H_E^\transpose\boldsymbol\lambda.
\end{equation}

For functions \(A,B\colon\cD_E\to\mathbb R\) satisfying \(A\leq B\), define
\begin{equation}\label{eq:L}
 L^{(A,B)}(R)
 =
 \sum_{\substack{\mathbf z\in\cD_E\\m_R(\mathbf z)>0}}
 m_R(\mathbf z)B(\mathbf z)
 +
 \sum_{\substack{\mathbf z\in\cD_E\\m_R(\mathbf z)<0}}
 m_R(\mathbf z)A(\mathbf z).
\end{equation}

\begin{proposition}[Support-function identity]
\label{prop:support}
Let \(R\in\cR_E\), and let
\[
 \mathbf a=(A(\mathbf z_s))_{s=1}^N,
 \qquad
 \mathbf b=(B(\mathbf z_s))_{s=1}^N.
\]
Then
\begin{equation}\label{eq:support}
 L^{(A,B)}(R)
 =
 \max_{\mathbf a\leq\mathbf c\leq\mathbf b}
 \mathbf m_R^\transpose\mathbf c.
\end{equation}
\end{proposition}

\begin{proof}
For a positive coefficient the maximum is attained at the upper bound,
and for a negative coefficient it is attained at the lower bound.
\end{proof}

We use the following inequality form of Farkas' lemma; see
\cite[Section~7.3]{Schrijver}.
\begin{theorem}[Farkas' lemma]\label{thm:farkas}
Let \(G\in\mathbb R^{p\times q}\) and \(\mathbf h\in\mathbb R^p\).
Exactly one of the following systems is feasible:
\begin{enumerate}
\item \(G\mathbf x\leq\mathbf h\), with \(\mathbf x\in\mathbb R^q\);
\item \(\mathbf y\geq0\), \(G^\transpose\mathbf y=0\), and
\(\mathbf h^\transpose\mathbf y<0\), with \(\mathbf y\in\mathbb R^p\).
\end{enumerate}
\end{theorem}
\smallskip

Let \(A,B\colon E\to[0,1]\) satisfy \(A\leq B\), and suppose that
the corresponding interval-valued data are quasi-copula compatible.
Let \(\underline Q_A\) and \(\overline Q_B\) be the extremal interval
envelopes from \eqref{eq:interval-lower}--\eqref{eq:interval-upper}, and
set
\begin{equation}\label{eq:grid-envelope-restrictions}
 A_E:=\underline Q_A|_{\cD_E},
 \qquad
 B_E:=\overline Q_B|_{\cD_E}.
\end{equation}
Then \(A_E,B_E\colon\cD_E\to[0,1]\) are discrete quasi-copulas and
\(A_E\leq B_E\). Put \(M:=|\cI_E|\).
Let
\[
 \mathbf a_E
 =
 \bigl(A_E(\mathbf z_1),\ldots,A_E(\mathbf z_N)\bigr)^\transpose,
 \qquad
 \mathbf b_E
 =
 \bigl(B_E(\mathbf z_1),\ldots,B_E(\mathbf z_N)\bigr)^\transpose.
\]
For \(\boldsymbol\lambda\in\mathbb R_{\geq0}^M\), define
\begin{align}\label{eq:Phi}
\begin{split}
 \Phi_E(\boldsymbol\lambda)
 &:=
 \sum_{(H_E^\transpose\boldsymbol\lambda)_s>0}
 (H_E^\transpose\boldsymbol\lambda)_s
 B_E(\mathbf z_s)
 +
 \sum_{(H_E^\transpose\boldsymbol\lambda)_s<0}
 (H_E^\transpose\boldsymbol\lambda)_s
 A_E(\mathbf z_s)\\
 &=\sum_{(H_E^\transpose\boldsymbol\lambda)_s>0}
 (H_E^\transpose\boldsymbol\lambda)_s(\mathbf b_E)_s
 +
 \sum_{(H_E^\transpose\boldsymbol\lambda)_s<0}
 (H_E^\transpose\boldsymbol\lambda)_s(\mathbf a_E)_s.
\end{split}
\end{align}
If \(\boldsymbol\lambda\in\mathbb N_0^M\) 
and $R_{\boldsymbol\lambda}$ is as in \eqref{def:R}, then
\begin{equation}\label{eq:Phi-L-relation}
 \Phi_E(\boldsymbol\lambda)
 =
 L^{(A_E,B_E)}(R_{\boldsymbol\lambda}).
\end{equation}

For each \(S\subseteq\{1,\ldots,N\}\), define the sign polytope
\begin{equation}\label{eq:sign-polytope}
 \cP_S
 :=
 \left\{
  \boldsymbol\lambda\in\mathbb R_{\geq0}^M:
  \mathbf 1^\transpose\boldsymbol\lambda=1,\ 
  (H_E^\transpose\boldsymbol\lambda)_s\geq0
       \text{ for }s\in S,\ 
  (H_E^\transpose\boldsymbol\lambda)_s\leq0
       \text{ for }s\notin S
 \right\}.
\end{equation}
Each nonempty \(\cP_S\) is a polytope and hence has finitely many
vertices; see, for example, \cite{Ziegler}. Its set of vertices is
denoted by \(\ver(\cP_S)\).

The following theorem gives the promised criterion for the original
finite interpolation problem.

\begin{theorem}[Finite copula interpolation criterion]
\label{thm:discrete-sandwich}
Let \(A,B\colon E\to[0,1]\) be as above, and let
\[
 A_E=\underline Q_A|_{\cD_E},
 \qquad
 B_E=\overline Q_B|_{\cD_E}.
\]
Then the following statements are equivalent:
\begin{enumerate}
\item There exists a \(d\)-copula \(C\) on \([0,1]^d\) such that
\[
 A\leq C|_E\leq B.
\]

\item There exists a discrete \(d\)-copula \(C_E\) on \(\cD_E\) such
that
\[
 A_E\leq C_E\leq B_E.
\]

\item
For every \(R\in\cR_E\),
\[
 L^{(A_E,B_E)}(R)\geq0.
\]

\item
For every \(S\subseteq\{1,\ldots,N\}\) such that
\(\cP_S\neq\varnothing\), and every
\(\boldsymbol\lambda\in\ver(\cP_S)\), one has
\[
 \Phi_E(\boldsymbol\lambda)\geq0.
\]
\end{enumerate}

In particular, finite interval-valued copula interpolation can be
decided by checking finitely many inequalities. Exact interpolation is
the special case \(A=B\).
\end{theorem}

\begin{proof}
If \textup{(1)} holds, then \(C_E:=C|_{\cD_E}\) is a discrete
\(d\)-copula. The extremal-envelope property in
\cref{prop:interval-envelopes} gives \(A_E\leq C_E\leq B_E\), so
\textup{(2)} holds. Conversely, if \textup{(2)} holds, then
\cref{prop:discrete-copula-extension} provides a \(d\)-copula extension
\(C\) of \(C_E\). Because \(E\subseteq\cD_E\) and the envelope bounds
enforce the original intervals at the points of \(E\), we have
\(A\leq C|_E\leq B\). Thus \textup{(1)} and \textup{(2)} are equivalent.

Let
\[
 \mathbf c
 =
 \bigl(C_E(\mathbf z_1),\ldots,C_E(\mathbf z_N)\bigr)^\transpose
\]
be the node-value vector of \(C_E\). Since \(A_E\) and \(B_E\) are
discrete quasi-copulas, they agree at every node where the copula
boundary conditions prescribe a fixed value. Moreover, nonnegative
elementary-cell volumes imply the lower-order monotonicity and
Lipschitz conditions; see \cref{lem:k-increasing-implies-quasi-copula}
below. Hence \textup{(2)} is equivalent to the feasibility of
\begin{equation}\label{eq:node-feasibility-extensions}
 H_E\mathbf c\geq\mathbf 0_M,
 \qquad
 \mathbf a_E\leq\mathbf c\leq\mathbf b_E,
\end{equation}
where vector inequalities are understood componentwise.

We next prove that \textup{(2)} implies \textup{(3)}. Suppose that
\(\mathbf c\) satisfies
\eqref{eq:node-feasibility-extensions}, and let
$
 \displaystyle 
 R
 =
 \bigsqcup_{\mathbf i\in\cI_E}
 \lambda_{\mathbf i}B_{\mathbf i}
 \in\cR_E.
$
After fixing an ordering of \(\cI_E\), write
$
 \boldsymbol\lambda
 =
 (\lambda_{\mathbf i})_{\mathbf i\in\cI_E}
 \in\mathbb N_0^M.
$
The multiplicity vector of \(R\) is
\[
 \mathbf m_R=H_E^\transpose\boldsymbol\lambda.
\]
It follows that
\[
 \mathbf m_R^\transpose\mathbf c
 =
 \boldsymbol\lambda^\transpose H_E\mathbf c
 \geq0.
\]
Moreover, since
\(\mathbf a_E\leq\mathbf c\leq\mathbf b_E\), the support-function
identity in \cref{prop:support} yields
\[
 \mathbf m_R^\transpose\mathbf c
 \leq
 \max_{\mathbf a_E\leq\mathbf x\leq\mathbf b_E}
 \mathbf m_R^\transpose\mathbf x
 =
 L^{(A_E,B_E)}(R).
\]
Consequently, \(L^{(A_E,B_E)}(R)\geq0\)
for every \(R\in\cR_E\). This proves that \textup{(2)} implies
\textup{(3)}.

We prove the converse by contraposition. Suppose that
\eqref{eq:node-feasibility-extensions} is infeasible. Write this system
as
\[
 G\mathbf c\leq\mathbf h,
\]
where
\begin{equation}\label{eq:farkas-block-matrix-extensions}
 G
 =
 \begin{bmatrix}
  -H_E\\
  I_N\\
  -I_N
 \end{bmatrix}
 \in\mathbb R^{(M+2N)\times N},
 \qquad
 \mathbf h
 =
 \begin{bmatrix}
  \mathbf 0_M\\
  \mathbf b_E\\
  -\mathbf a_E
 \end{bmatrix}
 \in\mathbb R^{M+2N}.
\end{equation}
By Farkas' lemma, in the form stated in \cref{thm:farkas}, there exist
\(\boldsymbol\lambda\in\mathbb R_{\geq0}^M\) and
\(\boldsymbol\mu,\boldsymbol\nu\in\mathbb R_{\geq0}^N\)
such that
\begin{equation}\label{eq:farkas-equality-extensions}
 -H_E^\transpose\boldsymbol\lambda
 +\boldsymbol\mu-\boldsymbol\nu
 =\mathbf 0_N
\end{equation}
and
\begin{equation}\label{eq:farkas-strict-extensions}
 \boldsymbol\mu^\transpose\mathbf b_E
 -
 \boldsymbol\nu^\transpose\mathbf a_E
 <0.
\end{equation}

Set \(\mathbf m=H_E^\transpose\boldsymbol\lambda\).
By \eqref{eq:farkas-equality-extensions},
\[
 \mathbf m=\boldsymbol\mu-\boldsymbol\nu.
\]
For each \(s\in\{1,\ldots,N\}\), if \(m_s\geq0\), then
\(\mu_s=m_s+\nu_s\), and hence
\begin{align*}
 \mu_s B_E(\mathbf z_s)-\nu_s A_E(\mathbf z_s)
 &=
 m_sB_E(\mathbf z_s)
 +
 \nu_s\bigl(B_E(\mathbf z_s)-A_E(\mathbf z_s)\bigr)
 \geq
 m_sB_E(\mathbf z_s).
\end{align*}
If \(m_s<0\), then \(\nu_s=\mu_s-m_s\), and therefore
\begin{align*}
 \mu_s B_E(\mathbf z_s)-\nu_s A_E(\mathbf z_s)
 &=
 m_sA_E(\mathbf z_s)
 +
 \mu_s\bigl(B_E(\mathbf z_s)-A_E(\mathbf z_s)\bigr)
 \geq 
 m_sA_E(\mathbf z_s).
\end{align*}
Summing over \(s\) and using
\eqref{eq:farkas-strict-extensions}, we obtain
\begin{equation}\label{eq:negative-Phi-extensions}
 \Phi_E(\boldsymbol\lambda)<0.
\end{equation}

The function
\(\Phi_E\colon\mathbb R_{\geq0}^M\to\mathbb R\) is continuous and
piecewise linear. It is also positively homogeneous of degree one:
\begin{equation}\label{eq:Phi-positive-homogeneity}
 \Phi_E(t\boldsymbol\lambda)=t\Phi_E(\boldsymbol\lambda)
 \qquad
 \text{for all }\boldsymbol\lambda\in\mathbb R_{\geq0}^M
 \text{ and }t\geq0.
\end{equation}
Indeed,
\(H_E^\transpose(t\boldsymbol\lambda)
=tH_E^\transpose\boldsymbol\lambda\), and multiplication by
\(t\geq0\) preserves the coordinate signs. Hence every term in
\(\Phi_E\) is multiplied by \(t\).

Since the
inequality in \eqref{eq:negative-Phi-extensions} is strict, there exists
a nonnegative rational vector
\(\boldsymbol\lambda'\in\mathbb Q_{\geq0}^M\) sufficiently close to
\(\boldsymbol\lambda\) such that
$
 \Phi_E(\boldsymbol\lambda')<0.
$
Choose a positive integer \(r\) such that
$
 \widehat{\boldsymbol\lambda}
 :=
 r\boldsymbol\lambda'
 \in\mathbb N_0^M.
$
By positive homogeneity,
$
 \Phi_E(\widehat{\boldsymbol\lambda})
 =
 r\Phi_E(\boldsymbol\lambda')
 <0.
$
The vector \(\widehat{\boldsymbol\lambda}\) defines the formal
combination
$
 R_{\widehat{\boldsymbol\lambda}}
 =
 \bigsqcup_{\mathbf i\in\cI_E}
 \widehat{\lambda}_{\mathbf i}B_{\mathbf i}
 \in\cR_E.
$
By \eqref{eq:Phi-L-relation},
\[
 L^{(A_E,B_E)}(R_{\widehat{\boldsymbol\lambda}})
 =
 \Phi_E(\widehat{\boldsymbol\lambda})
 <0.
\]
Thus, if \textup{(2)} fails, then \textup{(3)} fails. This proves the
equivalence of \textup{(2)} and \textup{(3)}.

It remains to prove the equivalence of \textup{(3)} and \textup{(4)}.
We first observe that \textup{(3)} is equivalent to
\begin{equation}\label{eq:Phi-cone-extensions}
 \Phi_E(\boldsymbol\lambda)\geq0
 \qquad
 \text{for every }
 \boldsymbol\lambda\in\mathbb R_{\geq0}^M.
\end{equation}
Indeed, by \eqref{eq:Phi-L-relation}, condition \textup{(3)} gives this
inequality for every
\(\boldsymbol\lambda\in\mathbb N_0^M\). Positive homogeneity then
extends it to every nonnegative rational vector, and continuity extends
it to every vector in \(\mathbb R_{\geq0}^M\). The converse follows by
restricting \eqref{eq:Phi-cone-extensions} to
\(\mathbb N_0^M\) and applying \eqref{eq:Phi-L-relation}.

Since \(\Phi_E\) is positively homogeneous, it is enough to check its
nonnegativity on the standard simplex
\[
 \Delta_M
 =
 \left\{
  \boldsymbol\lambda\in\mathbb R_{\geq0}^M:
  \mathbf 1^\transpose\boldsymbol\lambda=1
 \right\}.
\]
Indeed, if
\(\boldsymbol\lambda\in\mathbb R_{\geq0}^M\setminus\{\mathbf 0_M\}\),
then
$
 \widehat{\boldsymbol\lambda}
 =
 \frac{\boldsymbol\lambda}
      {\mathbf 1^\transpose\boldsymbol\lambda}
 \in\Delta_M
$
and
$
 \Phi_E(\boldsymbol\lambda)
 =
 \bigl(\mathbf 1^\transpose\boldsymbol\lambda\bigr)
 \Phi_E(\widehat{\boldsymbol\lambda}).
$
Moreover, \(\Phi_E(\mathbf 0_M)=0\).

The sign polytopes
$\cP_S$,
$S\subseteq\{1,\ldots,N\},$
cover \(\Delta_M\). Indeed, for each
\(\boldsymbol\lambda\in\Delta_M\), one may take
\[
 S(\boldsymbol\lambda)
 =
 \left\{
  s\in\{1,\ldots,N\}:
  (H_E^\transpose\boldsymbol\lambda)_s\geq0
 \right\},
\]
in which case
\(\boldsymbol\lambda\in\cP_{S(\boldsymbol\lambda)}\).

Fix \(S\subseteq\{1,\ldots,N\}\) such that
\(\cP_S\neq\varnothing\). On \(\cP_S\), the function \(\Phi_E\) is
given by the linear expression
\begin{equation}\label{eq:Phi-linear-on-sign-polytope}
 \Phi_E(\boldsymbol\lambda)
 =
 \sum_{s\in S}
 (H_E^\transpose\boldsymbol\lambda)_sB_E(\mathbf z_s)
 +
 \sum_{s\notin S}
 (H_E^\transpose\boldsymbol\lambda)_sA_E(\mathbf z_s).
\end{equation}
If a coordinate of \(H_E^\transpose\boldsymbol\lambda\) vanishes, its
contribution is zero, so its assignment to either sign class is
immaterial.

Since \(\Phi_E\) is linear on \(\cP_S\), it is nonnegative on
\(\cP_S\) if and only if it is nonnegative at every vertex of
\(\cP_S\). Therefore,
\[
 \Phi_E(\boldsymbol\lambda)\geq0
 \quad\text{for every }\boldsymbol\lambda\in\cP_S
\quad
\Longleftrightarrow
\quad
 \Phi_E(\boldsymbol\lambda)\geq0
 \quad\text{for every }
 \boldsymbol\lambda\in\ver(\cP_S).
\]
Since the sign polytopes \(\cP_S\) cover \(\Delta_M\), condition
\eqref{eq:Phi-cone-extensions} is equivalent to condition
\textup{(4)}. This proves the equivalence of \textup{(3)} and
\textup{(4)}.

Finally, there are only finitely many subsets
\(S\subseteq\{1,\ldots,N\}\), and each nonempty polytope \(\cP_S\) has
finitely many vertices. Hence condition \textup{(4)} consists of
finitely many inequalities.
\end{proof}

\medskip

\begin{example}[The complexity of the sign-polytope reduction] \label{ex:bivariate-sign-polytopes}
We now describe the vertices occurring in the sign-polytope reduction
for the grid from \cref{ex:bivariate-incidence}.

Write
\[
 \boldsymbol\lambda
 =
 \begin{pmatrix}
  \lambda_{11}&
  \lambda_{12}&
  \lambda_{13}&
  \lambda_{21}&
  \lambda_{22}&
  \lambda_{23}&
  \lambda_{31}&
  \lambda_{32}&
  \lambda_{33}
 \end{pmatrix}^\transpose
 \in\Delta_9,
\]
where
\[
 \Delta_9
 =
 \left\{
  \boldsymbol\lambda\in\mathbb R_{\geq0}^9:
  \sum_{i=1}^3\sum_{j=1}^3\lambda_{ij}=1
 \right\}.
\]
The coordinates correspond, respectively, to the elementary cells
\[
 B_{11},B_{12},B_{13},
 B_{21},B_{22},B_{23},
 B_{31},B_{32},B_{33}.
\]

Using the incidence matrix from \cref{ex:bivariate-incidence}, we
obtain
\begin{equation}\label{eq:example-H9-transpose-lambda}
 H_E^\transpose\boldsymbol\lambda
 =
 \begin{pmatrix}
 \lambda_{11}\\
 -\lambda_{11}+\lambda_{12}\\
 -\lambda_{12}+\lambda_{13}\\
 -\lambda_{13}\\
 -\lambda_{11}+\lambda_{21}\\
 \lambda_{11}-\lambda_{12}-\lambda_{21}+\lambda_{22}\\
 \lambda_{12}-\lambda_{13}-\lambda_{22}+\lambda_{23}\\
 \lambda_{13}-\lambda_{23}\\
 -\lambda_{21}+\lambda_{31}\\
 \lambda_{21}-\lambda_{22}-\lambda_{31}+\lambda_{32}\\
 \lambda_{22}-\lambda_{23}-\lambda_{32}+\lambda_{33}\\
 \lambda_{23}-\lambda_{33}\\
 -\lambda_{31}\\
 \lambda_{31}-\lambda_{32}\\
 \lambda_{32}-\lambda_{33}\\
 \lambda_{33}
 \end{pmatrix}.
\end{equation}
Thus, the sign-polytopal subdivision of \(\Delta_9\) is induced by
the sixteen coordinate hyperplanes of
\(H_E^\transpose\boldsymbol\lambda\), together with the nine facets
\[
 \lambda_{ij}=0,
 \qquad i,j\in\{1,2,3\}.
\]

A direct exact enumeration\footnote{The computation used Wolfram
Mathematica \cite{Wolfram2026}. Verification code is available at
\url{https://github.com/ZalarA/Copula-Interpolation}.} gives
\begin{equation}\label{eq:number-sign-vertices-example}
 \left|
 \bigcup_{\substack{S\subseteq\{1,\ldots,16\}\\
                    \cP_S\neq\varnothing}}
 \ver(\cP_S)
 \right|
 =
 159,
\end{equation}
where repetitions are omitted.

Enumerate the multiplier coordinates as
\[
 (\lambda_1,\ldots,\lambda_9)
 =
 (\lambda_{11},\lambda_{12},\lambda_{13},
  \lambda_{21},\lambda_{22},\lambda_{23},
  \lambda_{31},\lambda_{32},\lambda_{33}),
\]
and let \(\mathbf e_1,\ldots,\mathbf e_9\) denote the standard basis of
\(\mathbb R^9\). For every nonempty \(I\subseteq\{1,\ldots,9\}\), put
\[
 \mathbf v_I:=\frac1{|I|}\sum_{r\in I}\mathbf e_r.
\]
Of the \(159\) distinct vertices, \(147\) are uniform-support vertices
of the form \(\mathbf v_I\). If \(\mathcal I_r\) denotes the collection
of their supports of size \(r\), the enumeration gives
\[
\begin{array}{c|rrrrrrrrr|c}
 r & 1&2&3&4&5&6&7&8&9 & \text{total}\\ \hline
 |\mathcal I_r| & 9&12&22&20&32&24&18&9&1 & 147.
\end{array}
\]
The complete support lists are included with the verification code
cited above. Write
\[
 \cV_{\mathrm{unif}}
 :=
 \left\{\mathbf v_I:
 I\in\bigcup_{r=1}^9\mathcal I_r\right\}.
\]
There are twelve additional vertices that are not uniform on their
supports. Four of them have six nonzero coordinates:
\[
\begin{aligned}
\mathbf w_1&=\frac17(0,1,0,1,2,1,0,1,1)^\transpose,\qquad
\mathbf w_2=\frac17(0,1,0,1,2,1,1,1,0)^\transpose,\\
\mathbf w_3&=\frac17(0,1,1,1,2,1,0,1,0)^\transpose,\qquad
\mathbf w_4=\frac17(1,1,0,1,2,1,0,1,0)^\transpose.
\end{aligned}
\]
The remaining eight exceptional vertices have eight nonzero
coordinates:
\[
\begin{aligned}
\mathbf w_5&=\frac19(0,1,1,1,2,1,1,1,1)^\transpose,
&\qquad
\mathbf w_6&=\frac19(1,1,0,1,2,1,1,1,1)^\transpose,\\
\mathbf w_7&=\frac19(1,1,1,1,0,1,1,1,2)^\transpose,
&\qquad
\mathbf w_8&=\frac19(1,1,1,1,0,1,2,1,1)^\transpose,\\
\mathbf w_9&=\frac19(1,1,1,1,2,1,0,1,1)^\transpose,
&\qquad
\mathbf w_{10}&=\frac19(1,1,1,1,2,1,1,1,0)^\transpose,\\
\mathbf w_{11}&=\frac19(1,1,2,1,0,1,1,1,1)^\transpose,
&\qquad
\mathbf w_{12}&=\frac19(2,1,1,1,0,1,1,1,1)^\transpose.
\end{aligned}
\]

Consequently, the complete set of distinct sign-polytope vertices is
\[
 \boxed{
 \cV_{\mathrm{sign}}
 =
 \bigcup_{\substack{S\subseteq\{1,\ldots,16\}\\
                    \cP_S\neq\varnothing}}
 \ver(\cP_S)
 =
 \cV_{\mathrm{unif}}
 \cup
 \{\mathbf w_1,\ldots,\mathbf w_{12}\}.
 }
\]
In particular,
$
 |\cV_{\mathrm{sign}}|
 =
 147+12
 =
 159.
$
\end{example}
\medskip

\begin{example}[An explicit test of the discrete copula sandwich criterion] \label{ex:bivariate-noncopula-bounds}
We exhibit lower and upper extremal quasi-copula envelopes such that neither envelope is a discrete copula, although a discrete copula exists between them.

Let
\[
 E
 =
 \left\{0,\frac14,\frac23,1\right\}
 \times
 \left\{0,\frac13,\frac34,1\right\}.
\]
Thus, \(E\) is the entire grid from
\cref{ex:bivariate-incidence}, and the grid generated by \(E\) is
$
 \cD_E=E.
$
In particular, the incidence matrix \(H_E\) and the set of \(159\)
sign-polytope vertices from
\cref{ex:bivariate-incidence,ex:bivariate-sign-polytopes} remain
unchanged.

Using the node ordering from \cref{ex:bivariate-incidence}, define
\(A,B\colon E\to[0,1]\) by their node-value vectors
\[
 \mathbf a
 =
 \bigl(A(\mathbf z_1),\ldots,A(\mathbf z_{16})\bigr)^\transpose
 =
 \begin{pmatrix}
 0&
 0&
 0&
 0&
 0&
 0&
 \frac16&
 \frac14&
 0&
 \frac13&
 \frac5{12}&
 \frac23&
 0&
 \frac13&
 \frac34&
 1
 \end{pmatrix}^{\transpose}
\]
and
\[
 \mathbf b
 =
 \bigl(B(\mathbf z_1),\ldots,B(\mathbf z_{16})\bigr)^\transpose
 =
 \begin{pmatrix}
 0&
 0&
 0&
 0&
 0&
 0&
 \frac14&
 \frac14&
 0&
 \frac13&
 \frac12&
 \frac23&
 0&
 \frac13&
 \frac34&
 1
 \end{pmatrix}^{\transpose}.
\]
Clearly,
\[
 A\leq B
 \qquad\text{on }E.
\]

Both \(A\) and \(B\) satisfy the discrete quasi-copula boundary
conditions. Moreover, their increments along every elementary
horizontal or vertical edge are nonnegative and do not exceed the
length of the corresponding edge. Therefore, \(A\) and \(B\) are
discrete quasi-copulas.

Since the data are prescribed at every node of the generated grid,
the extremal envelopes reproduce the prescribed functions on
\(\cD_E\). Consequently,
\[
 A_E
 =
 \underline Q_A|_{\cD_E}
 =
 A,
 \qquad
 B_E
 =
 \overline Q_B|_{\cD_E}
 =
 B.
\]
Thus, their node-value vectors are
\[
 \mathbf a_E=\mathbf a,
 \qquad
 \mathbf b_E=\mathbf b.
\]

Neither of them is a discrete copula. Indeed,
\[
 \bigl(H_E\mathbf a_E\bigr)^\transpose
 =
 \begin{pmatrix}
 0&
 \frac16&
 \frac1{12}&
 \frac13&
 -\frac1{12}&
 \frac16&
 0&
 \frac13&
 0
 \end{pmatrix}
\]
and
\[
 \bigl(H_E\mathbf b_E\bigr)^\transpose
 =
 \begin{pmatrix}
 0&
 \frac14&
 0&
 \frac13&
 -\frac1{12}&
 \frac16&
 0&
 \frac14&
 \frac1{12}
 \end{pmatrix}.
\]
In both cases, the volume of the central elementary cell
\[
 B_{22}
 =
 \left[\frac14,\frac23\right]
 \times
 \left[\frac13,\frac34\right]
\]
is negative:
\[
 V_{A_E}(B_{22})
 =
 V_{B_E}(B_{22})
 =
 -\frac1{12}.
\]
Thus, neither \(A_E\) nor \(B_E\) is \(2\)-increasing.

We now use the sign-polytope reduction from
\cref{ex:bivariate-sign-polytopes}. Retain the notation
\[
 \cV_{\mathrm{sign}}
 =
 \cV_{\mathrm{unif}}
 \cup
 \{\mathbf w_1,\ldots,\mathbf w_{12}\}
\]
for the complete set of its \(159\) distinct vertices. For every
\(\boldsymbol\lambda\in\cV_{\mathrm{sign}}\), evaluate
\[
 \Phi_E(\boldsymbol\lambda)
 =
 \sum_{\bigl(H_E^\transpose\boldsymbol\lambda\bigr)_s>0}
 \bigl(H_E^\transpose\boldsymbol\lambda\bigr)_s
 B_E(\mathbf z_s)
 +
 \sum_{\bigl(H_E^\transpose\boldsymbol\lambda\bigr)_s<0}
 \bigl(H_E^\transpose\boldsymbol\lambda\bigr)_s
 A_E(\mathbf z_s).
\]
An exact calculation gives
\[
 \min_{\boldsymbol\lambda\in\cV_{\mathrm{sign}}}
 \Phi_E(\boldsymbol\lambda)
 =
 0.
\]
The minimum is attained precisely at
$\mathbf e_1$,
$\mathbf e_5$,
$\mathbf e_7$.
In particular,
\[
 \Phi_E(\boldsymbol\lambda)\geq0
 \qquad
 \text{for every }
 \boldsymbol\lambda\in\cV_{\mathrm{sign}}.
\]
It follows from \cref{thm:discrete-sandwich} that there exists a
discrete copula \(C_E\) satisfying
\[
 A_E\leq C_E\leq B_E.
\]
Thus, the sign-polytope inequalities provide a finite dual certificate
for the feasibility of the sandwich problem.

We may also recover a copula between the bounds directly by solving a
linear feasibility problem. Let
\[
 \mathbf c
 =
 \bigl(
 c(\mathbf z_1),\ldots,c(\mathbf z_{16})
 \bigr)^\transpose
 \in\mathbb R^{16}
\]
be the unknown vector of grid-node values. Consider the linear program
\begin{equation}\label{eq:example-copula-feasibility-LP}
\begin{aligned}
 \min_{\mathbf c\in\mathbb R^{16}}\quad
 &0\\
 \text{subject to}\quad
 &H_E\mathbf c\geq\mathbf0_9,\\
 &\mathbf a_E\leq\mathbf c\leq\mathbf b_E.
\end{aligned}
\end{equation}
The objective function is identically zero, so
\eqref{eq:example-copula-feasibility-LP} is simply a linear feasibility
problem. The inequality
\[
 H_E\mathbf c\geq\mathbf0_9
\]
requires all nine elementary-cell volumes to be nonnegative. The box
constraints require the node values to lie between \(A_E\) and \(B_E\).
No separate copula boundary constraints are needed because \(A_E\) and
\(B_E\) agree at every boundary node and satisfy the copula boundary
conditions there. Thus,
\[
 \mathbf a_E\leq\mathbf c\leq\mathbf b_E
\]
already forces \(\mathbf c\) to satisfy those boundary conditions.

Solving \eqref{eq:example-copula-feasibility-LP} gives, for example,
the feasible solution
\[
 \mathbf c_E
 =
 \begin{pmatrix}
 0&
 0&
 0&
 0&
 0&
 0&
 \frac16&
 \frac14&
 0&
 \frac13&
 \frac12&
 \frac23&
 0&
 \frac13&
 \frac34&
 1
 \end{pmatrix}^{\transpose}.
\]
It satisfies
\[
 \mathbf a_E\leq\mathbf c_E\leq\mathbf b_E.
\]
Moreover,
\[
 \bigl(H_E\mathbf c_E\bigr)^\transpose
 =
 \begin{pmatrix}
 0&
 \frac16&
 \frac1{12}&
 \frac13&
 0&
 \frac1{12}&
 0&
 \frac14&
 \frac1{12}
 \end{pmatrix}
 \geq\mathbf0_9^\transpose.
\]
Hence all elementary-cell volumes are nonnegative. Their sum is
\[
 0+\frac16+\frac1{12}+\frac13
 +0+\frac1{12}+0+\frac14+\frac1{12}
 =1.
\]
Since \(\mathbf c_E\) also satisfies the copula boundary conditions, it
is the node-value vector of a discrete copula \(C_E\) satisfying
\[
 A_E\leq C_E\leq B_E.
\]

This example illustrates both sides of the primal--dual description:
the nonnegativity of \(\Phi_E\) at the \(159\) sign-polytope vertices
provides a finite dual certificate of feasibility, while
\(\mathbf c_E\) provides an explicit primal witness.
\end{example}

\section{A sandwich criterion for
\texorpdfstring{\(k\)-increasing}{k-increasing} discrete quasi-copulas}
\label{sec:k-increasing-farkas}

We extend the incidence-matrix argument from \cref{sec:farkas} to
\(k\)-increasing discrete quasi-copulas, where \(2\leq k\leq d\).
Elementary \(d\)-cells are replaced by elementary \(k\)-dimensional
coordinate boxes. The existence of a \(k\)-increasing multivariate
function between prescribed bounds was studied in
\cite{KlementEtAlCoherence} through coherence and avoidance of sure
loss. Here we specialize the problem to finite data, derive the
corresponding Farkas certificate, and obtain a finite sign-polytope
test. For \(k=d\), the construction reduces to
\cref{thm:discrete-sandwich}.

Recall that
\[
 E=\{\mathbf x_1,\ldots,\mathbf x_n\}\subseteq[0,1]^d,
 \qquad
 \mathbf x_j=(x_{1j},\ldots,x_{dj}),
\]
is the finite set on which the lower and upper bounds are prescribed.
For each \(\ell\in\{1,\ldots,d\}\), the corresponding coordinate set is
\[
 D_\ell(E)=\{0,1,x_{\ell1},\ldots,x_{\ell n}\}
          =\{t_{\ell,0},\ldots,t_{\ell,n_\ell}\},
\]
where repeated values are omitted and
\[
 0=t_{\ell,0}<t_{\ell,1}<\cdots<t_{\ell,n_\ell}=1.
\]
The rectangular grid generated by \(E\) is
\[
 \cD_E=D_1(E)\times\cdots\times D_d(E).
\]

For the remainder of this section, fix \(k\in\{2,\ldots,d\}\). The
incidence matrix below records the volumes of all elementary
\(k\)-dimensional coordinate boxes in \(\cD_E\).

For every \(J\subseteq\{1,\ldots,d\}\) with \(|J|=k\), define
\[
 \cI_{E,J}
 =
 \prod_{\ell\in J}\{1,\ldots,n_\ell\},
 \qquad
 \cQ_{E,J}
 =
 \prod_{\ell\notin J}D_\ell(E).
\]

An elementary \(k\)-box of \(\cD_E\) is a set of the form
\begin{equation}\label{eq:elementary-k-box}
 B_{J,\mathbf i,\mathbf q}
 =
 \prod_{\ell=1}^d I_\ell,
\end{equation}
where
\[
 I_\ell
 =
 \begin{cases}
 [t_{\ell,i_\ell-1},t_{\ell,i_\ell}],
     & \ell\in J,\\
 \{q_\ell\},
     & \ell\notin J,
 \end{cases}
\]
for some
\[
 \mathbf i\in\cI_{E,J},
 \qquad
 \mathbf q\in\cQ_{E,J}.
\]
Thus, the coordinates indexed by \(J\) are active, whereas the remaining
coordinates are fixed.

Let
\begin{equation}\label{eq:k-cell-index-set}
 \cI_E^{(k)}
 =
 \left\{
  (J,\mathbf i,\mathbf q):
  J\subseteq\{1,\ldots,d\},\
  |J|=k,\
  \mathbf i\in\cI_{E,J},\
  \mathbf q\in\cQ_{E,J}
 \right\},
\end{equation}
and put
\[
 M_k=|\cI_E^{(k)}|.
\]
For
\(\alpha=(J,\mathbf i,\mathbf q)\in\cI_E^{(k)}\), write
\[
 B_\alpha=B_{J,\mathbf i,\mathbf q}.
\]
Its set of vertices consists of the \(2^k\) nodes obtained by choosing,
for each \(\ell\in J\), one of the two endpoints
\(t_{\ell,i_\ell-1}\) and \(t_{\ell,i_\ell}\), while keeping
\(z_\ell=q_\ell\) for \(\ell\notin J\).

Define the vertex-multiplicity function of \(B_\alpha\) by
\begin{equation}\label{eq:k-multiplicity}
 m_\alpha(\mathbf z)
 =
 \begin{cases}
 (-1)^{\#\{\ell\in J:
       z_\ell=t_{\ell,i_\ell-1}\}},
     & \mathbf z\in\ver(B_\alpha),\\
 0,  & \mathbf z\notin\ver(B_\alpha).
 \end{cases}
\end{equation}
After fixing an ordering of \(\cI_E^{(k)}\), define the
\(k\)-cell--vertex incidence matrix by
\begin{equation}\label{eq:H-k}
 H_E^{(k)}
 =
 \bigl[m_\alpha(\mathbf z_s)\bigr]_{
   \alpha\in\cI_E^{(k)},\;1\leq s\leq N}
 \in\{-1,0,1\}^{M_k\times N}.
\end{equation}

For a node-value vector
\[
 \mathbf c
 =
 \bigl(c(\mathbf z_1),\ldots,c(\mathbf z_N)\bigr)^\transpose,
\]
the component
\[
 \bigl(H_E^{(k)}\mathbf c\bigr)_\alpha
\]
is the \(k\)-dimensional \(C\)-volume of \(B_\alpha\). Hence
\begin{equation}\label{eq:k-increasing-matrix-condition}
 H_E^{(k)}\mathbf c\geq\mathbf 0_{M_k}
\end{equation}
is equivalent to the nonnegativity of the \(C\)-volume of every
elementary \(k\)-box. By finite additivity, this is equivalent to
\(k\)-increasingness on \(\cD_E\).
\medskip

\begin{lemma}
\label{lem:k-increasing-implies-quasi-copula}
Let
\[
 \cD=D_1\times\cdots\times D_d\subseteq[0,1]^d
\]
be a finite rectangular grid such that \(0,1\in D_\ell\) for every
\(\ell\in[d]\). Let \(F\colon\cD\to[0,1]\) be grounded and have
uniform one-dimensional margins. If \(F\) is \(k\)-increasing for some
\(k\in\{2,\ldots,d\}\), then \(F\) is \(r\)-increasing for every
\(r\in\{1,\ldots,k\}\) and, in particular, \(F\) is a discrete
quasi-copula.
\end{lemma}

\begin{proof}
We first show that \(r\)-increasingness implies
\((r-1)\)-increasingness for every \(r\in\{2,\ldots,d\}\).

Let \(B_J\) be an \((r-1)\)-dimensional coordinate box with
\(|J|=r-1\). If one of its fixed coordinates is zero, then every
vertex of \(B_J\) lies in a grounded face. Hence
$
 V_F(B_J)=0.
$

Suppose, therefore, that all fixed coordinates of \(B_J\) are
positive. Choose \(s\notin J\), and denote the fixed \(s\)-th
coordinate of \(B_J\) by \(q_s>0\). Let \(\widetilde B\) be the
\(r\)-dimensional coordinate box obtained from \(B_J\) by replacing
the fixed coordinate \(\{q_s\}\) with the interval \([0,q_s]\). By
the definition of the alternating volume,
\[
 V_F(\widetilde B)
 =
 V_F(B_J)-V_F(B_J^{(s=0)}),
\]
where \(B_J^{(s=0)}\) denotes the copy of \(B_J\) contained in the
face whose \(s\)-th coordinate is zero. Groundedness gives
$
 V_F(B_J^{(s=0)})=0.
$
Consequently,
\[
 V_F(B_J)=V_F(\widetilde B)\geq0.
\]
Thus \(r\)-increasingness implies \((r-1)\)-increasingness. Iterating
this argument shows that \(F\) is \(r\)-increasing for every
\(r\in\{1,\ldots,k\}\).

In particular, \(F\) is \(1\)-increasing and hence increasing in each
variable. It remains to verify the discrete Lipschitz condition. Let
\(\mathbf u,\mathbf v\in\cD\) differ only in coordinate \(i\), with
\(u_i\leq v_i\). Since \(F\) is increasing,
\[
 F(\mathbf v)-F(\mathbf u)\geq0.
\]

Moreover, \(F\) is \(2\)-increasing. For every \(j\neq i\),
\(2\)-increasingness implies that the increment of \(F\) over
\([u_i,v_i]\) is nondecreasing in the \(j\)-th coordinate. Indeed, if
\(u_j\leq u_j'\), then applying \(2\)-increasingness to the coordinate
box with active coordinates \(i,j\), intervals
$[u_i,v_i]$,
$[u_j,u_j']$,
and all remaining coordinates fixed gives
\[
\begin{aligned}
&F(u_1,\ldots,v_i,\ldots,u_j,\ldots,u_d)
 -
 F(u_1,\ldots,u_i,\ldots,u_j,\ldots,u_d)\\
&\qquad\leq
 F(u_1,\ldots,v_i,\ldots,u_j',\ldots,u_d)
 -
 F(u_1,\ldots,u_i,\ldots,u_j',\ldots,u_d).
\end{aligned}
\]
Raising all coordinates other than \(i\) successively to \(1\), we
obtain
\[
\begin{aligned}
 F(\mathbf v)-F(\mathbf u)
 &\leq
 F(1,\ldots,1,v_i,1,\ldots,1)
 -
 F(1,\ldots,1,u_i,1,\ldots,1).
\end{aligned}
\]
By the uniform one-dimensional marginal condition, the right-hand
side is \(v_i-u_i\). Therefore,
\[
 0\leq F(\mathbf v)-F(\mathbf u)\leq v_i-u_i.
\]
Thus \(F\) satisfies the discrete quasi-copula edge conditions. Hence
\[
 \lvert F(\mathbf u)-F(\mathbf v)\rvert
 \leq
 \sum_{\ell=1}^d\lvert u_\ell-v_\ell\rvert
\]
for all \(\mathbf u,\mathbf v\in\cD\), and \(F\) is a discrete
quasi-copula.
\end{proof}

\begin{remark}
\label{rem:one-increasing-not-sufficient}
The assumption \(k\geq2\) in
\cref{lem:k-increasing-implies-quasi-copula} is essential.
For \(k=1\), groundedness, uniform one-dimensional margins, and
\(1\)-increasingness imply that the function is coordinatewise
nondecreasing, but they do not imply the Lipschitz condition.
Consequently, a grounded and standardized \(1\)-increasing function
need not be a discrete quasi-copula.

For example, let
\[
 D=\left\{0,\frac25,\frac12,1\right\},
 \qquad
 \cD=D^2,
\]
and define \(F\colon\cD\to[0,1]\) by the usual groundedness and
marginal conditions,
\[
 F(x,0)=F(0,y)=0,
 \qquad
 F(x,1)=x,
 \qquad
 F(1,y)=y,
\]
and by
\[
\begin{array}{c|cc}
 F(x,y) & x=\frac25 & x=\frac12\\ \hline
 y=\frac25 & 0 & \frac25\\[1mm]
 y=\frac12 & 0 & \frac25
\end{array}
\]
at the remaining interior nodes. Then \(F\) is grounded, has uniform
one-dimensional margins, and is increasing in each coordinate.
However,
\[
 F\left(\frac12,\frac25\right)
 -
 F\left(\frac25,\frac25\right)
 =
 \frac25
 >
 \frac1{10}
 =
 \frac12-\frac25.
\]
Thus \(F\) fails the discrete Lipschitz condition and is not a
discrete quasi-copula. Hence \(1\)-increasingness alone is not
sufficient; the implication in
\cref{lem:k-increasing-implies-quasi-copula} relies on the resulting
\(2\)-increasingness when \(k\geq2\).
\end{remark}

\medskip

We next introduce formal nonnegative integer combinations of elementary
\(k\)-boxes:
\begin{equation}\label{eq:formal-family-k}
 \cR_E^{(k)}
 =
 \left\{
  \bigsqcup_{\alpha\in\cI_E^{(k)}}
  \lambda_\alpha B_\alpha:
  \boldsymbol\lambda
  =
  (\lambda_\alpha)_{\alpha\in\cI_E^{(k)}}
  \in\mathbb N_0^{M_k}
 \right\}.
\end{equation}
As before, \(\bigsqcup\) denotes a formal multiset union.

For
\[
 R
 =
 \bigsqcup_{\alpha\in\cI_E^{(k)}}
 \lambda_\alpha B_\alpha
 \in\cR_E^{(k)},
\]
define
\begin{equation}\label{eq:formal-k-multiplicity}
 m_R^{(k)}(\mathbf z)
 =
 \sum_{\alpha\in\cI_E^{(k)}}
 \lambda_\alpha m_\alpha(\mathbf z),
 \qquad
 \mathbf z\in\cD_E.
\end{equation}
The corresponding multiplicity vector is
\begin{equation}\label{eq:k-multiplicity-vector}
 \mathbf m_R^{(k)}
 =
 \bigl(
  m_R^{(k)}(\mathbf z_1),\ldots,
  m_R^{(k)}(\mathbf z_N)
 \bigr)^\transpose
 =
 \bigl(H_E^{(k)}\bigr)^\transpose\boldsymbol\lambda.
\end{equation}

For functions \(A,B\colon\cD_E\to\mathbb R\) satisfying \(A\leq B\),
define
\begin{equation}\label{eq:L-k}
 L_k^{(A,B)}(R)
 =
 \sum_{\substack{\mathbf z\in\cD_E\\
                  m_R^{(k)}(\mathbf z)>0}}
 m_R^{(k)}(\mathbf z)B(\mathbf z)
 +
 \sum_{\substack{\mathbf z\in\cD_E\\
                  m_R^{(k)}(\mathbf z)<0}}
 m_R^{(k)}(\mathbf z)A(\mathbf z).
\end{equation}

The proof of the following identity is identical to that of
\cref{prop:support}.

\begin{proposition}[\(k\)-dimensional support-function identity]
\label{prop:support-k}
Let \(R\in\cR_E^{(k)}\), and let
\[
 \mathbf a
 =
 \bigl(A(\mathbf z_1),\ldots,A(\mathbf z_N)\bigr)^\transpose,
 \qquad
 \mathbf b
 =
 \bigl(B(\mathbf z_1),\ldots,B(\mathbf z_N)\bigr)^\transpose.
\]
Then
\begin{equation}\label{eq:support-k}
 L_k^{(A,B)}(R)
 =
 \max_{\mathbf a\leq\mathbf c\leq\mathbf b}
 \bigl(\mathbf m_R^{(k)}\bigr)^\transpose\mathbf c.
\end{equation}
\end{proposition}

\begin{proof}
The maximum is attained coordinatewise: at \(b_s\) when
\(m_R^{(k)}(\mathbf z_s)>0\), and at \(a_s\) when
\(m_R^{(k)}(\mathbf z_s)<0\). This is the same argument as in the proof
of \cref{prop:support}.
\end{proof}

Let \(A,B\colon E\to[0,1]\) satisfy \(A\leq B\), and suppose that
the corresponding interval-valued data are quasi-copula compatible.
Let \(\underline Q_A\) and \(\overline Q_B\) be the extremal envelopes
from \eqref{eq:interval-lower}--\eqref{eq:interval-upper}, and set
\begin{equation}\label{eq:k-grid-envelope-restrictions}
 A_E:=\underline Q_A|_{\cD_E},
 \qquad
 B_E:=\overline Q_B|_{\cD_E}.
\end{equation}
Then \(A_E\) and \(B_E\) are discrete quasi-copulas satisfying
\(A_E\leq B_E\). The pointwise argument from
\cref{prop:interval-envelopes}, which uses only monotonicity and the
Lipschitz inequalities at the relevant points, also applies on the
grid. Hence every discrete quasi-copula \(Q_E\) on \(\cD_E\) satisfies
\begin{equation}\label{eq:k-extremal-extension-property}
 A\leq Q_E|_E\leq B
 \quad\Longleftrightarrow\quad
 A_E\leq Q_E\leq B_E.
\end{equation}
In particular, this equivalence holds for every \(k\)-increasing
discrete quasi-copula on \(\cD_E\).

Set
\[
 \mathbf a_E
 =
 \bigl(
  A_E(\mathbf z_1),\ldots,
  A_E(\mathbf z_N)
 \bigr)^\transpose,
 \qquad
 \mathbf b_E
 =
 \bigl(
  B_E(\mathbf z_1),\ldots,
  B_E(\mathbf z_N)
 \bigr)^\transpose.
\]

For
\(\boldsymbol\lambda\in\mathbb R_{\geq0}^{M_k}\), define
\begin{align}
\label{eq:Phi-k}
 \Phi_E^{(k)}(\boldsymbol\lambda)
 &:=
 \sum_{\bigl((H_E^{(k)})^\transpose
       \boldsymbol\lambda\bigr)_s>0}
 \bigl((H_E^{(k)})^\transpose
       \boldsymbol\lambda\bigr)_s
 B_E(\mathbf z_s)
 +
 \sum_{\bigl((H_E^{(k)})^\transpose
       \boldsymbol\lambda\bigr)_s<0}
 \bigl((H_E^{(k)})^\transpose
       \boldsymbol\lambda\bigr)_s
 A_E(\mathbf z_s).
\end{align}
If
\(\boldsymbol\lambda\in\mathbb N_0^{M_k}\), then
\begin{equation}\label{eq:Phi-L-k-relation}
 \Phi_E^{(k)}(\boldsymbol\lambda)
 =
 L_k^{(A_E,B_E)}
 \bigl(R_{\boldsymbol\lambda}\bigr),
\end{equation}
where
$
 R_{\boldsymbol\lambda}
 =
 \bigsqcup_{\alpha\in\cI_E^{(k)}}
 \lambda_\alpha B_\alpha.
$

For each \(S\subseteq\{1,\ldots,N\}\), define the sign polytope
\begin{equation}\label{eq:k-sign-polytope}
 \cP_S^{(k)}
 =
 \left\{
  \boldsymbol\lambda\in\mathbb R_{\geq0}^{M_k}:
  \mathbf 1^\transpose\boldsymbol\lambda=1,\
  \bigl((H_E^{(k)})^\transpose\boldsymbol\lambda\bigr)_s\geq0
       \text{ for }s\in S,\
  \bigl((H_E^{(k)})^\transpose\boldsymbol\lambda\bigr)_s\leq0
       \text{ for }s\notin S
 \right\}.
\end{equation}

\begin{theorem}[\(k\)-increasing discrete sandwich criterion]
\label{thm:k-increasing-discrete-sandwich}
Fix \(k\in\{2,\ldots,d\}\). Let \(A,B\colon E\to[0,1]\) be as above,
and let
\[
 A_E=\underline Q_A|_{\cD_E},
 \qquad
 B_E=\overline Q_B|_{\cD_E}.
\]
Then the following statements are equivalent:
\begin{enumerate}
\item
There exists a \(k\)-increasing discrete quasi-copula \(C_E\) on
\(\cD_E\) such that
\[
 A_E\leq C_E\leq B_E.
\]
Equivalently, there exists such a quasi-copula whose restriction to
\(E\) satisfies \(A\leq C_E|_E\leq B\).

\item
For every \(R\in\cR_E^{(k)}\),
\[
 L_k^{(A_E,B_E)}(R)\geq0.
\]

\item
For every \(S\subseteq\{1,\ldots,N\}\) such that
\(\cP_S^{(k)}\neq\varnothing\), and every
\(\boldsymbol\lambda\in\ver(\cP_S^{(k)})\), one has
\[
 \Phi_E^{(k)}(\boldsymbol\lambda)\geq0.
\]
\end{enumerate}

In particular, the existence of a \(k\)-increasing discrete
quasi-copula between \(A_E\) and \(B_E\) can be decided by checking
finitely many inequalities.
\end{theorem}

\begin{proof}
The proof follows the proof of \cref{thm:discrete-sandwich}, with
\[
 H_E,\quad
 M,\quad
 \cR_E,\quad
 \mathbf m_R,\quad
 L^{(A_E,B_E)},\quad
 \Phi_E,\quad
 \cP_S
\]
replaced, respectively, by
\[
 H_E^{(k)},\quad
 M_k,\quad
 \cR_E^{(k)},\quad
 \mathbf m_R^{(k)},\quad
 L_k^{(A_E,B_E)},\quad
 \Phi_E^{(k)},\quad
 \cP_S^{(k)}.
\]

Indeed, by \eqref{eq:k-extremal-extension-property}, assertion
\textup{(1)} is equivalent to the existence of a \(k\)-increasing
discrete quasi-copula \(C_E\) on \(\cD_E\) satisfying
$
 A_E\leq C_E\leq B_E.
$
Let
\[
 \mathbf c
 =
 \bigl(C_E(\mathbf z_1),\ldots,C_E(\mathbf z_N)\bigr)^\transpose
\]
be its node-value vector. Since \(A_E\) and \(B_E\) are discrete
quasi-copulas, they agree at all nodes at which the quasi-copula
boundary conditions prescribe a fixed value. Consequently, the
inequalities
\[
 \mathbf a_E\leq\mathbf c\leq\mathbf b_E
\]
force the corresponding function \(C_E\) to be grounded and to have
uniform one-dimensional margins.

Moreover,
\[
 H_E^{(k)}\mathbf c\geq\mathbf0_{M_k}
\]
is equivalent to \(k\)-increasingness on \(\cD_E\), because its
components are the volumes of the elementary \(k\)-boxes and every
\(k\)-dimensional coordinate box with vertices in \(\cD_E\) is a
finite union of elementary \(k\)-boxes. Since \(k\geq2\),
\cref{lem:k-increasing-implies-quasi-copula} implies that the
corresponding function \(C_E\) is a discrete quasi-copula.

Thus, assertion \textup{(1)} is equivalent to the feasibility of 
\begin{equation}\label{eq:k-node-feasibility} 
    H_E^{(k)}\mathbf c\geq\mathbf 0_{M_k}, 
    \qquad 
    \mathbf a_E \leq\mathbf c \leq\mathbf b_E. 
\end{equation}

If \(\mathbf c\) satisfies \eqref{eq:k-node-feasibility}, then for every
\(R\in\cR_E^{(k)}\), with coefficient vector
\(\boldsymbol\lambda\in\mathbb N_0^{M_k}\),
\[
 \bigl(\mathbf m_R^{(k)}\bigr)^\transpose\mathbf c
 =
 \boldsymbol\lambda^\transpose
 H_E^{(k)}\mathbf c
 \geq0.
\]
Together with \cref{prop:support-k}, this gives
\[
 L_k^{(A_E,B_E)}(R)\geq0.
\]

Conversely, if \eqref{eq:k-node-feasibility} is infeasible, apply
Farkas' lemma in the form stated in \cref{thm:farkas} to
\[
 G_k
 =
 \begin{bmatrix}
  -H_E^{(k)}\\
  I_N\\
  -I_N
 \end{bmatrix},
 \qquad
 \mathbf h_k
 =
 \begin{bmatrix}
  \mathbf 0_{M_k}\\
  \mathbf b_E\\
  -\mathbf a_E
 \end{bmatrix}.
\]
Exactly as in the proof of \cref{thm:discrete-sandwich}, the resulting
Farkas multiplier yields
\[
 \boldsymbol\lambda\in\mathbb R_{\geq0}^{M_k}
 \qquad\text{such that}\qquad
 \Phi_E^{(k)}(\boldsymbol\lambda)<0.
\]
The function \(\Phi_E^{(k)}\) is continuous, piecewise linear, and
positively homogeneous of degree one. Approximating
\(\boldsymbol\lambda\) by a nonnegative rational vector and clearing
denominators therefore produces
\(\widehat{\boldsymbol\lambda}\in\mathbb N_0^{M_k}\) such that
\[
 L_k^{(A_E,B_E)}
 \bigl(R_{\widehat{\boldsymbol\lambda}}\bigr)
 =
 \Phi_E^{(k)}(\widehat{\boldsymbol\lambda})
 <0.
\]
This proves the equivalence of \textup{(1)} and \textup{(2)}.

Finally, positive homogeneity reduces the nonnegativity of
\(\Phi_E^{(k)}\) on \(\mathbb R_{\geq0}^{M_k}\) to its nonnegativity on
the standard simplex
\[
 \left\{
  \boldsymbol\lambda\in\mathbb R_{\geq0}^{M_k}:
  \mathbf 1^\transpose\boldsymbol\lambda=1
 \right\}.
\]
The sign polytopes \(\cP_S^{(k)}\) cover this simplex, and
\(\Phi_E^{(k)}\) is linear on each of them. It is therefore
nonnegative on every \(\cP_S^{(k)}\) if and only if it is nonnegative
at every vertex of every nonempty \(\cP_S^{(k)}\). This proves the
equivalence of \textup{(2)} and \textup{(3)}.

Since there are finitely many sign polytopes and each has finitely many
vertices, only finitely many inequalities need to be checked.
\end{proof}

\begin{remark}
    For \(k=d\), the only active coordinate set is
    \(J=\{1,\ldots,d\}\). Hence
    $H_E^{(d)}=H_E$,
    $\cR_E^{(d)}=\cR_E$,
    and \cref{thm:k-increasing-discrete-sandwich} reduces to
    \cref{thm:discrete-sandwich}.
\end{remark}

\section{Conclusion and outlook}
\label{sec:concluding-remarks}

We established a finite criterion for exact and interval-valued copula
interpolation on an arbitrary finite set \(E\subseteq[0,1]^d\). The
first step is an explicit pair of extremal quasi-copula envelopes. The
coordinates of \(E\) then generate a canonical rectangular grid on
which the interpolation problem becomes
\[
 H_E\mathbf c\geq\mathbf0,
 \qquad
 \mathbf a_E\leq\mathbf c\leq\mathbf b_E.
\]
Farkas' lemma identifies the \(L^{(A_E,B_E)}\)-criterion with a dual
infeasibility certificate for this system. Positive homogeneity and
piecewise linearity further reduce the infinite family of formal-cell
conditions to a finite test at sign-polytope vertices.

The same construction applies to \(k\)-increasing discrete
quasi-copulas after replacing elementary \(d\)-cells by elementary
\(k\)-dimensional coordinate boxes. For \(k=d\), it recovers the copula
criterion.

The reduction is finite but may still be large: both the canonical grid
and the collection of sign-polytope vertices can grow rapidly with the
dimension and the number of distinct data coordinates. A natural next
step is therefore to compare the node-value formulation with the
alternative cell-mass formulation, exploit redundancy and symmetry in
the incidence matrices, and develop algorithms that generate only the
dual certificates relevant to a given data set.

\end{document}